\documentclass[preprint]{imsart}

\RequirePackage[OT1]{fontenc}
\RequirePackage{amsthm,amsmath}
\RequirePackage{natbib}
\RequirePackage[colorlinks,citecolor=blue,urlcolor=blue]{hyperref}

\arxiv{arXiv:0000.0000}  

\startlocaldefs
\numberwithin{equation}{section}
\theoremstyle{plain}

\endlocaldefs

\usepackage{mathpazo}
\usepackage{courier}
\usepackage{amsmath,amsfonts,amssymb,amsthm}
\usepackage{graphicx}
\usepackage{natbib,url}
\usepackage{color}
\usepackage{float}

\newtheorem{rmq}{Remark}[section]

\DeclareMathOperator{\clr}{CLR}

\DeclareMathOperator{\crit}{crit}

\DeclareMathOperator{\logclr}{logCLR}
\newcommand{\LCLR}{logCLR}

\usepackage[draft]{fixme}
\fxusetheme{color}
\FXRegisterAuthor{h}{ah}{H}
\FXRegisterAuthor{p}{ap}{P}
\FXRegisterAuthor{a}{aa}{A}

\begin{document}
\begin{frontmatter}
\title{
Clustering transformed compositional data using $K$-means, with applications in gene expression and bicycle sharing system data 
}

\runtitle{
Clustering transformed compositional data using $K$-means}

\begin{aug}
\author{\fnms{Antoine} \snm{Godichon-Baggioni}\thanksref{t1}\ead[label=e1]{godichon@insa-toulouse.fr}},
\author{\fnms{Cathy} \snm{Maugis-Rabusseau}\ead[label=e2]{cathy.maugis@insa-toulouse.fr}}
\and
\author{\fnms{Andrea} \snm{Rau}\ead[label=e3]{andrea.rau@inra.fr}}

\thankstext{t1}{Corresponding author}
\runauthor{A. Godichon-Baggioni et al.}

\affiliation{\thanksmark{m1} Institut de Math{\'e}matiques de Toulouse, Universit{\'e} Toulouse III - Paul Sabatier}
\address{ Institut de Math{\'e}matiques de Toulouse\\ 
 Université Toulouse III - Paul Sabatier\\
118 route de Narbonne\\ 
31062 Toulouse, France.\\
\printead{e1}\\}

\affiliation{\thanksmark{m2} Institut de Math\'ematiques de Toulouse, INSA de Toulouse}
\address{Institut Math\'ematiques de Toulouse\\
 INSA de Toulouse,\\ 
 135 avenue de Rangueil,\\ 
 31077 Toulouse, France.\\
\printead{e2}\\}

\affiliation{\thanksmark{m3} GABI, INRA, AgroParisTech, Universit{\'e} Paris-Saclay}
\address{Institut National de la Recherche Agronomique\\ 
Domaine de Vilvert\\
78352 Jouy-en-Josas, France.\\
\printead{e3}\\}
\end{aug}

\begin{abstract}
Although there is no shortage of clustering algorithms proposed in the literature, the question of the most relevant strategy for clustering compositional data (i.e., data made up of profiles, whose rows belong to the simplex) remains largely unexplored in cases where the observed value of an observation is equal or close to zero for one or more samples. This work is motivated by the analysis of two sets of compositional data, both focused on the categorization of profiles but arising from considerably different applications: (1) identifying groups of co-expressed genes from high-throughput RNA sequencing data, in which a given gene may be completely silent in one or more experimental conditions; and (2) finding patterns in the usage of stations over the course of one week in the Velib' bicycle sharing system in Paris, France. For both of these applications, we focus on the use of appropriately chosen data transformations, including the Centered Log Ratio and a novel extension we propose called the Log Centered Log Ratio, in conjunction with the $K$-means algorithm. We use a nonasymptotic penalized criterion, whose penalty is calibrated with the slope heuristics, to select the number of clusters present in the data. Finally, we illustrate the performance of this clustering strategy, which is implemented in the Bioconductor package \texttt{coseq}, on both the gene expression and bicycle sharing system data.
\end{abstract}



\end{frontmatter}

\section{Introduction}

Compositional data are made up of the relative proportions of a whole and can be represented in the simplex of $d$ parts:
\begin{equation*}
\mathcal{S}^{d}:= \left\lbrace x=\left( x_{1},...,x_{d} \right) \in \mathbb{R}^{d}|  \sum_{i=1}^{d}x_{i} = 1 , x_{i}> 0, \forall i \right\rbrace .
\end{equation*}
Such data are a common phenomenon in several research domains, including evolutionary ecology \citep{aebischer1993compositional,bingham2007misclassified}, geochemistry \citep{buccianti2006compositional,miesch1977log}, economics \citep{desarbo1995analyzing,longford2006stability}, and genomic surveys \citep{Friedman2012}; see \cite{pawlowsky2011compositional} for more examples of applications. The statistical analysis of compositional data has been the focus of research for well over a century \citep{pearson1896mathematical,chayes1960correlation}. One of the initial concerns was to propose appropriate methods able to account for the nature of compositional data, notably the fact that they are subject to a unit sum constraint, which invalidates many standard statistical approaches. For example, \cite{pawlowsky2001geometric} provided a formal definition of metric center and variance for random compositional data, as well as a law of large numbers. \cite{mateu2013normal} provided a definition for the normal distribution in the simplex, as well as the relevant central limit theorem. The concept of proximity between two compositional vectors has been defined using several proposed distance measures, such as Aitchison's distance \citep{aitchison1982statistical}, which is typically based on data transformations such as the Centered Log Ratio (CLR), the Additive Log Ratio (ALR), or the Isometric Log Ratio (ILR) \citep{egozcue2003isometric}. These compositional data transformations  facilitate the use of statistical methods based on Gaussian distributions and Euclidean structures for compositional data.

Many clustering methods exist in the literature, and they can primarily be divided into two classes: model-based methods, such as mixture models \citep{mclachlan2004finite}, and methods based on dissimilarity distances, such as hierarchical clustering \citep{ward1963hierarchical}, $K$-means \citep{macqueen1967some}, or $K$-medians \citep{cardot2012fast}. However, despite the large number of existing clustering methods, to our knowledge there has been relatively little attention paid to the most appropriate strategy for clustering  compositional data \citep{tauber1999spurious, zhou1991logratio, martin1998critical}. In recent work, \cite{rau2017transformation} proposed the use of data transformations (either the arcsine or logit) and Gaussian mixture models to cluster compositional data arising from high-throughput transcriptome sequencing data. This strategy led to satisfactory results in practice as the proposed transformations removed the linear dependence present among data coordinates, thus enabling estimation of model parameters; however, such an approach implies that the dependencies among coordinates are entirely ignored. 

In this work, our aim is twofold: (1) to introduce clustering methods for compositional data based on dissimilarity distances, in particular the $K$-means algorithm \citep{macqueen1967some}; and (2) to use or define appropriate transformations for compositional data that account for the dependencies among coordinates. In particular, we investigate three clustering strategies for the task of clustering compositional data via the $K$-means algorithm with the usual Euclidean distance: (1) using untransformed data, which is arguably the most intuitive but does not directly account for the compositional nature of the data; (2) using data transformed using the CLR \citep{aitchison1982statistical}, which is specific to data on the simplex; and (3) transforming data using a novel extension of the CLR transformation, called the Log Centered Log Ratio (\LCLR), consisting of a modification designed to provide greater separation for edge case observations, i.e. those with compositional values close to zero in one or more coordinates and those located near the vertices of the simplex. Finally, we select the number of clusters using a nonasymptotic criterion whose penalty is calibrated using the slope heuristics \citep{baudry2012slope, fischer2011number}. This approach of using a $K$-means algorithm in conjunction with compositional data transformations is implemented in the Bioconductor package \texttt{coseq}, and we apply it to two practical problems: identifying groups of co-expressed genes from high-throughput RNA sequencing data, and finding patterns in the usage of stations in a bicycle sharing system.  

The remainder of this paper is organized as follows: in Section \ref{pres}, we describe the context and data from the transcriptomic and bicycle sharing applications in greater detail. In Section \ref{methods}, we present the use of the $K$-means algorithm for the three proposed transformations, as well as the model selection approach used to select the number of clusters.  A full analysis using this approach is performed on the transcriptomic and bicycle sharing system data in Sections \ref{rnaseq} and \ref{bss}, respectively. Finally, we provide some conclusions and recommendations in Section \ref{conclusion}.

\section{Description of motivating data}\label{pres}
\subsection{Transcriptomic data}\label{presrna}

Gene expression studies now routinely make use of high-throughput sequencing technology, which directly sequences reverse-transcribed RNA molecules in an approach called RNA sequencing (RNA-seq). After aligning sequenced reads to a reference genome and quantifying  the number of reads attributed to each gene, RNA-seq data correspond to tables of read counts or pseudocounts $\left( Y_{i,j} \right)$ representing the number of sequenced reads observed for genes $i=1,\ldots,n$ in biological samples $j=1,\ldots,d$, the latter of which may arise from several experimental conditions (e.g., across time, in different tissues). In this work, our focus will be on identifying groups of {\it co-expressed} genes with the same expression dynamics across all biological samples from an RNA-seq study; these co-expression clusters are often assumed to be involved in similar biological processes or to be candidates for co-regulation.

RNA-seq data tend to be characterized by highly skewed count values covering several orders of magnitude, as well as variable {\it library sizes} among samples (i.e., the total number of sequenced reads in a sample). It is often assumed that the read count $Y_{i,j}$ is proportional to the expression of gene $i$, weighted by the library size of sample $j$ (as gene read counts in samples with larger total numbers of sequenced reads tend to have larger read counts) as well as its length (i.e., the number of nucleotides making up its coding region, as longer genes also tend to have larger read counts). As such, clustering strategies for RNA-seq data must account for both the length of each gene and the library size of each sample. Regarding the former, as is typically done in the context of differential analyses, we will make use of per-sample scaling normalization factors $\left(t_{j}\right)$ calculated using the Trimmed Mean of M-values (TMM) approach \citep{Robinson2010}. For all $j=1,\ldots,d$, let
\begin{align*}
& \ell_{j}:= t_{j}\sum_{i=1}^{n}Y_{i,j},
& s_{j}:= \frac{\ell_{j}}{\sum_{j'=1}^{d}\ell_{j'}/d},
\end{align*}
be the normalized library sizes and the associated scaling factors by which raw counts are divided, respectively. As such, read counts normalized for differences in library sizes may be calculated as $Y_{i,j} / s_{j}$.
To account for differences in the length of each gene, similarly to \cite{rau2017transformation}, we calculate the normalized expression \textit{profiles} defined for each gene $i$ by $X_{i}=\left( X_{i,1},\ldots,X_{i,d}\right)$, where
\begin{equation*}
X_{i,j} := \frac{Y_{i,j}/s_{j}+1}{\sum_{j'=1}^{d}\left( Y_{i,j'}/s_{j'}+1\right)},\quad \quad \forall j = 1,\ldots,d,
\end{equation*}
and a constant of 1 has been added to the normalized expression $Y_{i,j} / s_{j}$ prior to calculating the profiles due to the presence of 0's in the data. Note that we now have a table of compositional values $\left( X_{i,j} \right)$ whose rows $X_{i}$ belong to the simplex $\mathcal{S}^{d}$. In other words, the normalized expression profile value for a given gene in each sample is now relative to the total number of reads observed across all samples, meaning that this measure is independent of both its absolute expression level and its length; the normalized expression profiles thus facilitate the clustering of expression dynamics, as desired. 

In Section \ref{rnaseq}, we will focus our attention on two sets of RNA-seq data:
\begin{itemize}
\item {\bf Embryonic mouse neocortex:} \cite{fietz2012transcriptomes} studied the expansion of the neocortex in five embryonic (day 14.5) mice by analyzing the transcriptome of the ventricular zone (VZ), subventricular zone (SVZ) and cortical place (CP). Details on the sample preparation, sequencing, and quantification are provided in the Supplementary Materials of \citet{fietz2012transcriptomes}. In the current work, raw read counts were downloaded from the Digital Expression Explorer \citep{Ziemann2015} as described in \cite{rau2017transformation}. The data consist of transcriptome-wide measurements in three tissues with five replicates each. As suggested by \cite{rau2017transformation}, clustering will be performed on the full set of data (after filtering, on $n=8969$ genes), and visualization of resulting clusters will be done on profiles summed over each tissue.

\item \textbf{Dynamic expression in embryonic flies:} \cite{graveley2011developmental} characterized the expression dynamics of the fly using RNA-seq over 27 stages of development, ranging from early embryo to adult males and females. Raw read counts were obtained using the online ReCount resource \citep{Frazee2011}. We focus here on a subset of these data arising from 12 embryonic samples collected at 2-hour intervals over a 24-hour period, with a single replicate for each time point. After filtering, we obtain a subset of $n=9524$ genes.
\end{itemize}

\begin{figure}[th!]
\includegraphics[scale=0.8]{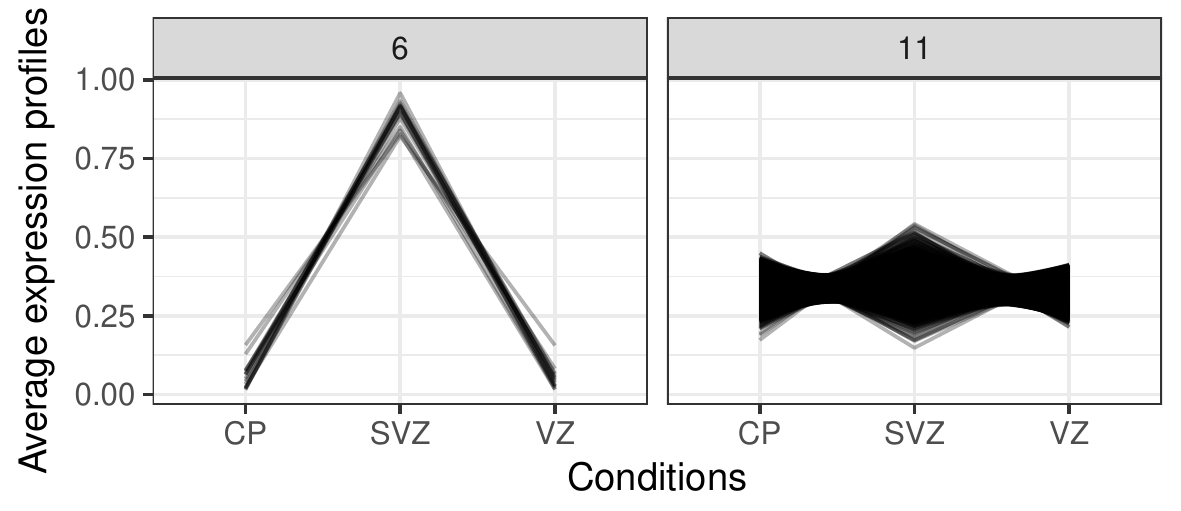}
\caption{Example of two groups of genes from the mouse neocortex RNA-seq data displaying normalized expression profiles that are highly-specific to the second tissue (left) and largely unchanged across tissues (right).
\label{ex_prof}
}
\includegraphics[scale=0.8]{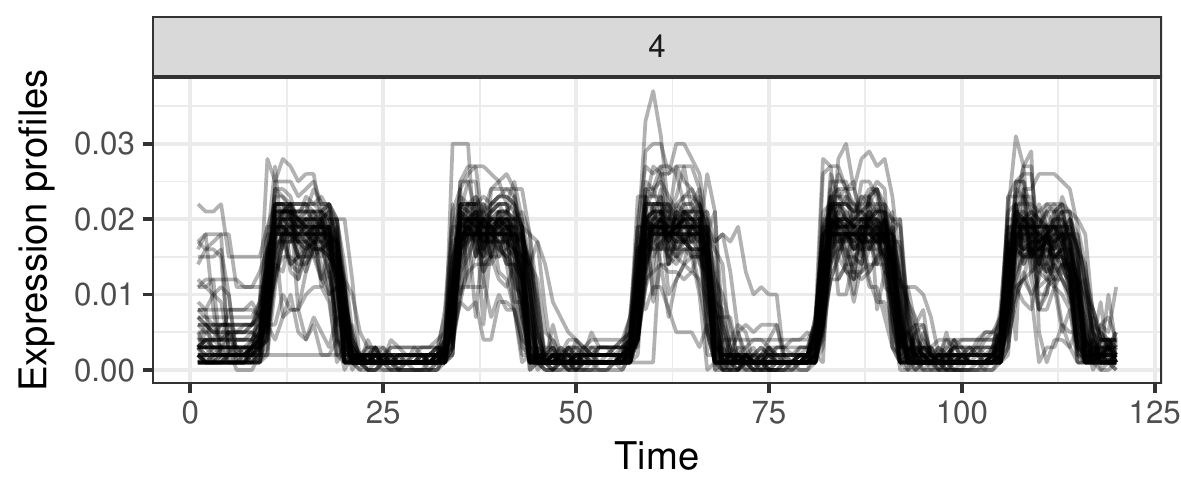}
\caption{Example of hourly bicycle station occupancy profiles (collected over a five-week time period) displaying a strong periodic trend in the Velib' bicycle sharing system data. \label{ex_prof_velib}}
\end{figure}

For both of the RNA-seq datasets described above, there are a large number of genes whose expression remains unchanged in different tissues or at different time points (i.e., whose profiles are close to the center of the simplex), and a relatively small number of genes with expression highly specific to a single tissue or time point (i.e., whose profiles are close to a vertex of the simplex); see Figure \ref{ex_prof} for an example from the mouse neocortex data. In the remainder of the paper, we refer to the latter as observations with {\it highly-specific} profiles, and we will focus on proposing a clustering method able to highlight these small but important groups of genes.

\subsection{Bicycle sharing data}
The Velib' bicycle sharing system was introduced in Paris, France in 2007 and is made up of around 14,500 bicycles in 1230 rental stations located around the Paris metropolitan area. Velib' subscribers can check out and return bicycles from any rental stations, and a fleet of vehicules provides a daily overnight redistribution of bicycles among rental stations. Bicycle rental stations located above 60 meters of elevation belong to a special category called V+ stations.

Hourly Velib' station occupancy data, in terms of available bicycles and docks, were downloaded from a period covering 1 week (between 11am on Sunday, August 31st and 11pm on Sunday, September 7th, 2014) from open-data APIs provided by the JCDecaux company, as described in \cite{bouveyron2015discriminative}. The data correspond to count tables $\left( Y_{i,j}\right)$, representing the raw occupancy counts (i.e., the number of available bicycles) for stations $i=1,\ldots,n$ at times $j=1,\ldots,d$ over the course of this week. Due to the daily habits of Velib' users (e.g., going to work), we note that the station occupancy counts tend to display periodic patterns across time points (see Figure \ref{ex_prof_velib}). In recent work, to account for differences in the size of rental stations (i.e. the number of docking points),  \cite{bouveyron2015discriminative} proposed dividing each row of counts by the total capacity of the station to yield occupancy proportions; they then clustered the rental stations in terms of the number of available docks at any given time. In our work, rather than using the fullness of each station, we instead focus on the dynamics of users' habits across time (e.g., whether bicycles tend to be rented and returned during working hours).  In a similar manner to the RNA-seq data described above, we thus calculate the relative occupancy of each Velib' station across the full time period studied; this facilitates clustering stations according to specific user behaviors across time while also accounting for differences in capacity among stations. We thus calculate the occupancy profiles $\left( X_{i} \right)_{i}$, defined for all $i$ by $X_{i}= \left( X_{i,1},\ldots,X_{i,d}\right)$ where
\begin{equation*}
X_{i,j}:= \frac{Y_{i,j}+1}{ \sum_{j'=1}^{d} \left( Y_{i,j'}+1\right)}, \quad \forall j=1,\ldots,d.
\end{equation*}
Note that as with the RNA-seq data, a constant of 1 is added to occupancy counts due to the presence of 0's. 
We thus obtain a proportion table $\left( X_{i,j} \right)$, whose rows $X_{i}$ belong to the simplex $\mathcal{S}^{d}$. 

In Section \ref{bss}, we focus on the occupancy profiles from the Velib' bicycle sharing data. These data are available in the R package \texttt{funFEM} \citep{bouveyron2015discriminative}, and correspond to counts of available bicycles in $n=1213$ stations at each hour, during one week ($d=183$). As weekend rental habits tend to be much different from those during weekdays, we focus in particular on the occupancy counts from Monday $12$am to Friday $11$pm ($d=120$). Although clustering will be performed on the full set of data, results are visualized after summing across days for each time interval ($d=24$) since a strong periodicity in station occupancy was observed from day to day.

\section{Methods}\label{methods}
In what follows, we consider a table of compositional data $X= \left(X_{i,j} \right)$, with $i=1,\ldots,n$ and $j=1,\ldots,d$. Each row $X_{i}$ of this table, which we also refer to as \textit{profiles}, belongs to the simplex $\mathcal{S}^{d}$.

\subsection{$K$-means algorithm}
Let $X_{1},\ldots,X_{n}$ be a set of $d$-dimensional points to be clustered into $K$ clusters. We consider the usual Euclidean norm $\| . \|_{2}$. Let $\mathcal{C}^{(K)}= \left\lbrace C_{k},k=1,\ldots,K \right\rbrace$ be a partition into $K$ clusters, and let $\mu_{k}$ be the mean of the cluster $C_{k}$: 
\[
\mu_{k}:= \frac{1}{\vert C_{k} \vert}\sum_{i \in C_{k}}X_{i},
\]
where $\vert C_{k} \vert$ is the cardinality of cluster $k$. 
The aim of $K$-means is to minimize the sum of squared errors (SSE), defined for each set of clusters $\mathcal{C}^{(K)}$ by
\begin{equation*}
\text{SSE} \left( \mathcal{C}^{(K)}\right) := \sum_{k=1}^{K}\sum_{i \in C_{k}} \left\| X_{i} - \mu_{k} \right\|_{2}^{2} ,
\end{equation*}
with $i \in C_{k}$ if $\left\| X_{i} - \mu_{k} \right\|_{2} = \min_{k'=1,\ldots,K} \left\| X_{i}- \mu_{k'} \right\|_{2}$. Many algorithms have been proposed to implement $K$-means clustering, and we consider the well-known one introduced by \cite{macqueen1967some}. Note that minimizing the SSE is known to be a NP-hard problem, and as such the $K$-means algorithm can only converge to a local minimum. In what follows, we will consider the $K$-means based on the Euclidean distance, using either untransformed or transformed data.  

\medskip

\noindent\textbf{Transformations for compositional data:}
Let $h: \mathcal{S}^{d} \longrightarrow \mathbb{R}^{d-1}$ be a differentiable and bijective function. The mapping $h$ enables the definition of a Euclidean structure on the simplex \citep{pawlowsky2001geometric}. As previously noted, the most usual transformations $h$ in the literature for compositional data are the CLR, the ALR and the ILR \citep{aitchison1982statistical}. As these different transformations yielded similar results for the two applications considered here, we focus on the $\clr :\mathcal{S}^{d} \longrightarrow \mathbb{R}^{d}$, defined for all $x \in \mathcal{S}^{d}$ by
\begin{align}
\clr(x) & := \left( \ln \left( \frac{x_{1}}{g(x)} \right),\ldots, \ln \left( \frac{x_{d}}{g(x)} \right) \right),  \label{eqn:clrdef}
\end{align}
where $g(x)$ is the geometric mean of $x$. In this case, the transformed values do not belong to $\mathbb{R}^{d-1}$ but to the hyperplane of $\mathbb{R}^{d}$ with normal vector $\left( 1,\ldots,1\right)$. Note that there is a strict equivalence between clustering untransformed compositional data with the $K$-means algorithm using Aitchison's distance \citep{aitchison1982statistical}, and clustering CLR-transformed data with the $K$-means algorithm using the Euclidean distance, i.e minimizing
\[
\text{SSE}_{\clr}\left( \mathcal{C}^{(K)} \right) := \sum_{k=1}^{K}\sum_{i\in C_{k}} \left\| \clr \left( X_{i} \right) - \mu_{k,\clr} \right\|_{2}^{2},
\]
where $\mu_{k,\clr} $ is the arithmetic mean of the CLR-transformed data belonging to  cluster $C_{K}$:
\[
\mu_{k,\clr} := \frac{1}{\vert C_{k} \vert}\sum_{i \in C_{k}}\clr \left( X_{i} \right) .
\]
Remark that $\clr^{-1}\left( \mu_{k,\clr} \right)$ corresponds to the center of the cluster when Aitchison's distance is used \citep{pawlowsky2001geometric}.

\bigskip

\begin{figure}[th!]
\centering
\includegraphics[height=7.5cm]{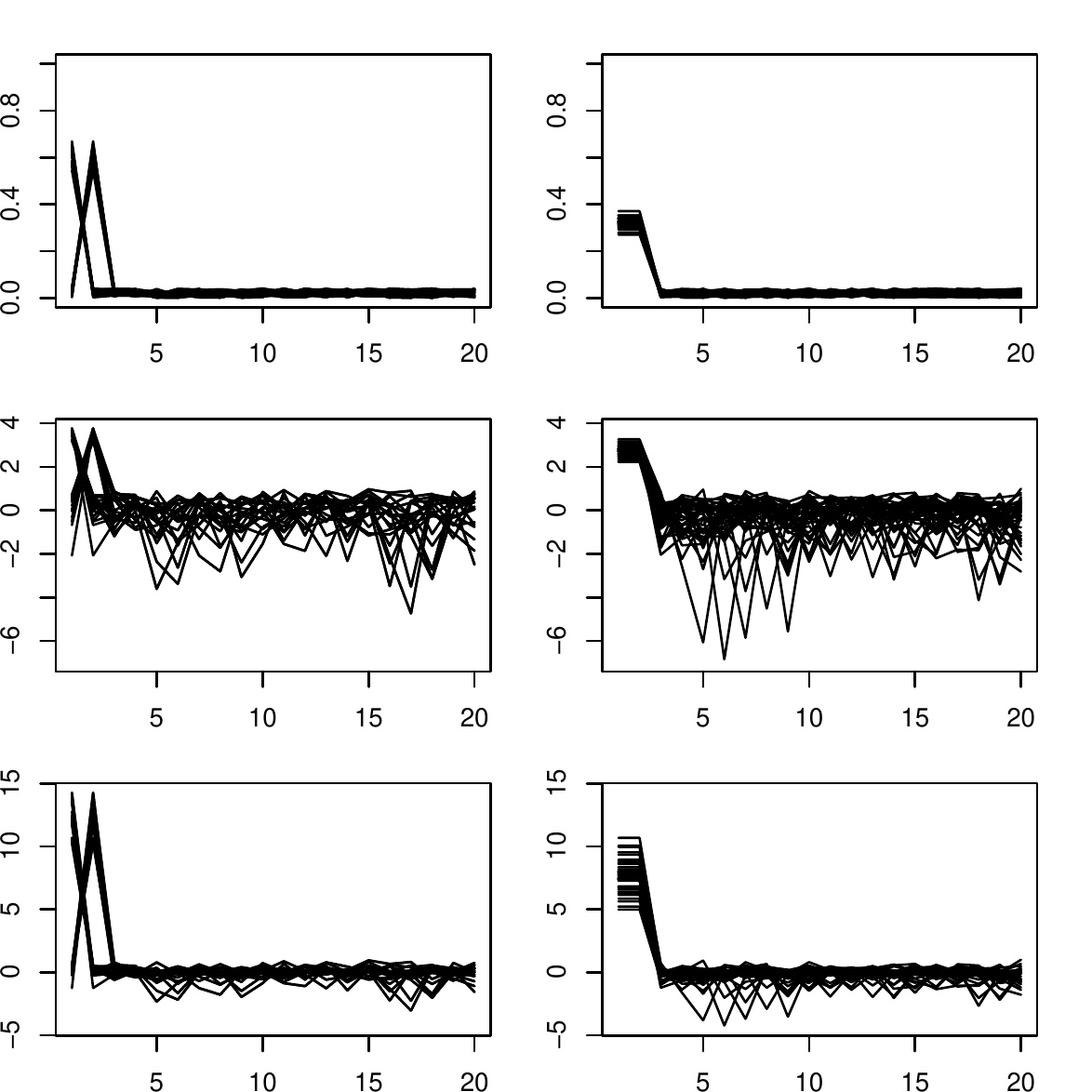}
\caption{Simulated highly-specific profiles on the simplex, on untransformed data (top) and after transformation using the CLR (middle) or \LCLR\ (bottom). On the left, two clusters of profiles that are specific to a single sample (either the first or second) are included; on the right, a single cluster of profiles specific to the first pair of samples is included. \label{figlogclr}}
\end{figure}

\noindent\textbf{Log Centered Log Ratio transformation:}
For data of moderate dimension, when a large number of coordinates have very small proportions, the CLR transformation tends to be quite sensitive to small fluctuations close to zero (see Figure \ref{figlogclr} for an example on simulated data). This can have a strong undesired effect on clustering results (see Section \ref{rnaseq}) when a small number of observations have highly-specific profiles. To account for this phenomenon by giving more importance to coordinates with large relative values, we propose a novel extension of the CLR for compositional data called the Log Centered Log Ratio (\LCLR). For all  $x \in \mathcal{S}^{d}$, the \LCLR\  is defined by $\logclr (x) := \left( \logclr\left( x_{1}\right),\ldots,\logclr\left( x_{d}\right) \right)$, where for all $j$,
$$
\logclr \left( x_{j}\right) := \left\{
    \begin{array}{ll}
        - \left[ \ln \left( 1- \ln \left[ x_{j}/g(x) \right] \right) \right]^{2}  & \mbox{if } x_{j}/g(x) \leq 1, \\
        \left( \ln \left[ x_{j}/g(x) \right] \right)^{2} & \mbox{otherwise,}
    \end{array}
\right.
$$
and $g(x)$ is the geometric mean of $x$. The additional $\log$ term when $\frac{x_{j}}{g(x)} \leq 1$ accords less importance in the transformation to samples with relatively weak proportions, while the squared term facilitates the concentration of profiles close to the center of the simplex $\left( \frac{1}{d},\ldots,\frac{1}{d}\right)$ (see Figure \ref{figlogclr}). Performing $K$-means clustering with this transformation amounts to minimizing
\[
\text{SSE}_{\logclr}\left( \mathcal{C}^{(K)}  \right) := \sum_{k=1}^{K}\sum_{i\in C_{k}} \left\| \logclr \left( X_{i} \right) - \mu_{k,\logclr} \right\|_{2}^{2},
\]
where $\mu_{k,\logclr} $ is the arithmetic mean of the transformed data belonging to the cluster $C_{K}$:
\[
\mu_{k,\logclr} := \frac{1}{\vert C_{k} \vert}\sum_{i \in C_{k}}\logclr\left( X_{i}\right) .
\]  

\begin{rmq}
Note that for all positive constants $p$, it is also possible to replace the exponent $2$ in the definition of the \LCLR\ transformation by an exponent of $p$.
\end{rmq}

\subsection{Choice of the number of clusters}\label{secrit}
Many criteria exist in the literature to select the number of clusters for $K$-means \citep{calinski1974dendrite,krzanowski1988criterion}, and the most common are the Gap statistic \citep{tibshirani2001estimating} and the maximization of the averaged Silhouette \citep{kaufman2009finding}. However, in our case, these methods tend to lead to the choice of an overly parsimonious model in practice. We instead focus on the selection method introduced by \cite{fischer2011number}, which is based on the minimization of a nonasymptotic penalized criterion defined for all positive integers $K$:
\begin{align*}
\crit (K) :& = \sum_{k=1}^{K}\sum_{i \in C_{k}}\left\| h \left( X_{i} \right) - \mu_{k,h} \right\|_{2}^{2} + \text{pen}\left( K \right) \\
& = \sum_{i=1}^{n}\min_{k=1,\ldots,K} \left\| h \left( X_{i} \right) - \mu_{k,h} \right\|_{2}^{2} + \text{pen}(K),
\end{align*}
where $h$ is the transformation used (identity for untransformed data). The penalty function $\text{pen}: \mathbb{N} \longrightarrow \mathbb{R}_{+}$ is defined by
\begin{equation*}
\text{pen}(K) := a_{h} \sqrt{Knd }.
\end{equation*}
The penalty term is known up to a multiplicative constant $a_{h}$, which in practice is calibrated using the slope heuristics method \citep{birge2007minimal,baudry2012slope} implemented in the R package \texttt{capushe} \citep{brault2011package}. The number of selected clusters is thus
\[
\widehat{K}:= \arg \min_{K\leq n} \crit (K).
\]

\subsection{Conditional probabilities}
By associating a clustering partition $\mathcal{C}^{(K)}$ obtained with the $K$-means algorithm with results obtained with the EM algorithm \citep{DLR77} for spherical Gaussian mixture models, the conditional probability that each observation $X_{i}$ belongs to a cluster $C_{k}$ may be calculated as follows:
\begin{equation*}
\tau_{i,k,h}:=\frac{ \left| C_{k} \right| \phi \left(h \left( X_{i} \right); \mu_{k,h},\sigma_{k,h}^2I_{d}\right)}{\sum_{k'=1}^{K}\left| C_{k'} \right| \phi \left(h \left( X_{i} \right);  \mu_{k',h},\sigma_{k',h}^2I_{d} \right)},
\end{equation*}
where $\phi \left(.; \mu_{k,h},\sigma_{k,h}^{2}I_{d}\right)$ is the density function of a Gaussian vector of mean $\mu_{k,h}$ and variance $\sigma_{k,h}^{2}I_{d}$, and $\sigma_{k,h}^{2}$ is the within-cluster variance:
\[
\sigma_{k,h}^{2} = \frac{1}{\vert C_{k} \vert}\sum_{i \in C_{k}}\left\| h \left( X_{i} \right) - \mu_{k,h} \right\|_{2}^{2}.
\] 
Note that alternative model selection criteria, such as the Bayesian Information Criterion (BIC) or Integrated Completed Likelihood (ICL) criterion \citep{schw_1978, Biernacki2000}, could also be used to select the number of clusters. Nevertheless, as for the use of the EM algorithm for Gaussian mixture models, for the studied data, the BIC and ICL criteria often lead to the choice of a very large number of clusters for the $K$-means algorithm. 

\subsection{Implementation in \texttt{coseq}}
We have implemented a user-friendly interface for the $K$-means clustering approach described above for the CLR and \LCLR\ transformations in the R/Bioconductor package \texttt{coseq}. Model selection is provided via the slope heuristics as described above, and a variety of customizable graphics of cluster profiles may be easily generated using dedicated plotting functions. Example scripts illustrating clustering analyses for profiles have been included directly in the package vignette.

\section{RNA-seq data results}\label{rnaseq}

We applied the $K$-means algorithm using the Euclidean distance described in the previous section to the two RNA-seq datasets described in Section \ref{pres} (referred to as the mouse or fly data, respectively), using either untransformed, CLR-transformed, or \LCLR-transformed profiles for $K=2,\ldots,40$ clusters. For each dataset and each transformation,  the nonasymptotic penalized criterion described in Section~\ref{secrit} is used to identify the appropriate number of clusters $\widehat{K}$ (see Table \ref{nbcluster}); diagnostic plots for the calibration of the slope heuristics penalty as well as profile plots for all identified clusters under each transformation are shown in Appendix. We note that alternative model selection criteria, such as the Gap statistic and averaged Silhouette, led to models with very small numbers of clusters (2 or 3). As such, the following results correspond to the models selected via the nonasymptotic penalized criterion.

\begin{table}[th!]\centering
\begin{tabular}{l|ccc}
& Identity  & CLR & \LCLR \\
\hline
Mouse & 22  & 20& 20 \\
Fly & 20  & 17 & 21
\end{tabular}
\caption{Number of clusters $\widehat{K}$ selected for each RNA-seq dataset and each proposed transformation, where model selection was performed by minimizing the penalized criterion defined in Section \ref{secrit} and calibrated with the slope heuristics. \label{nbcluster}}
\end{table}

\medskip

\begin{figure}[th!]\centering
\includegraphics[height=8.25cm]{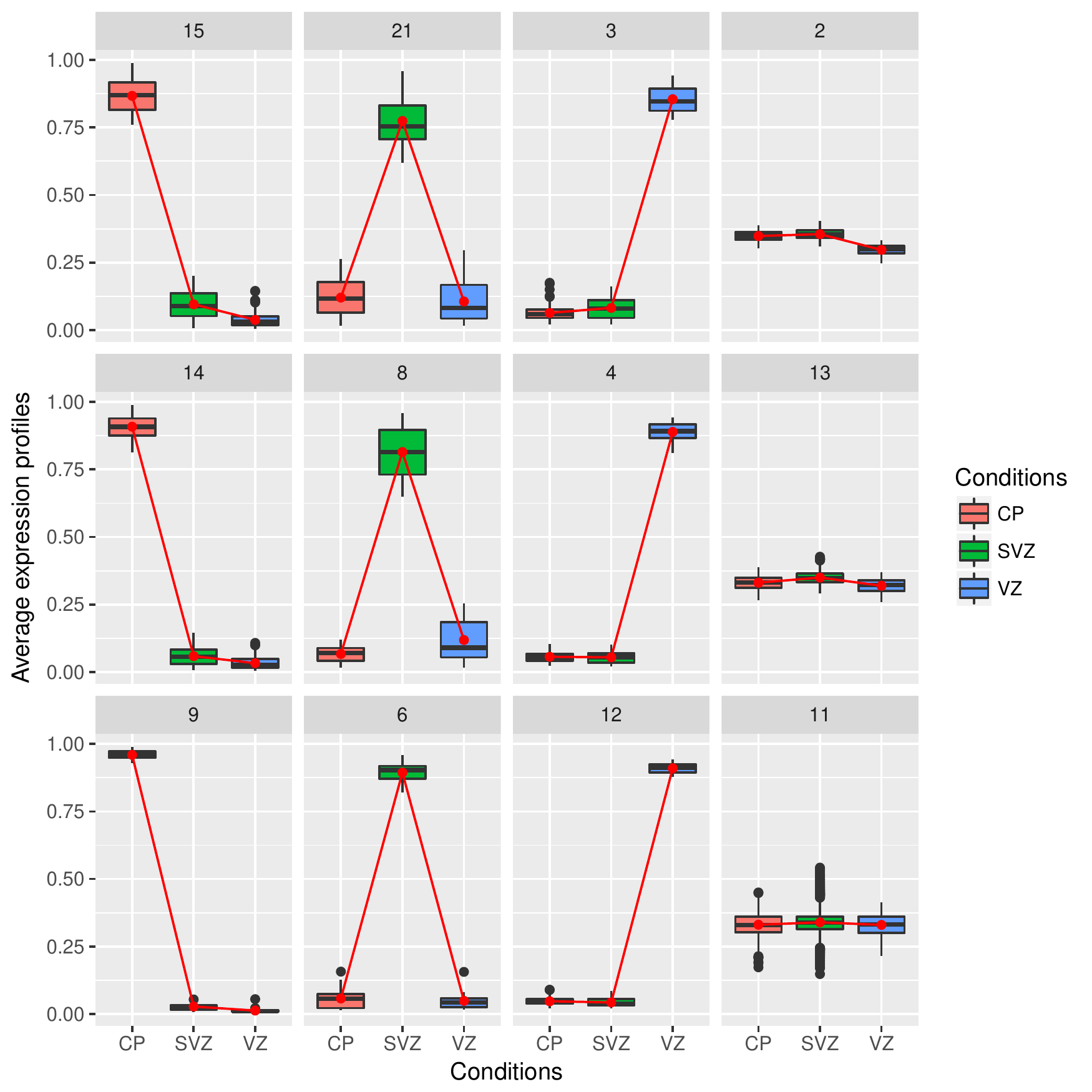}\caption{Examples of clusters of genes from the mouse RNA-seq data: per-cluster normalized expression profiles of selected clusters obtained with the $K$-means algorithm with no transformation (top), and with the CLR (middle) and \LCLR\ (bottom) transformations are shown. Analogous clusters for each approach have been vertically aligned. Connected red lines correspond to the mean profiles for each condition. \label{mouse}}
\end{figure}

By examining the clusters of profiles identified in each dataset under each transformation, we may note qualitative differences in the results for each approach, in particular for genes with highly-specific profiles. In the case of the mouse RNA-seq data, comparisons among analagous profile clusters for untransformed, CLR-, and \LCLR-transformed data indicate that untransformed and CLR-transformed data (Figure \ref{mouse}, top and middle) tend to yield a larger number of small, less-variable clusters situated at the center of the simplex, corresponding to several distinct clusters representing genes that are largely nondifferentially expressed. These two strategies also tend to yield more diffuse clusters for profiles close to the vertices of the simplex (Figure \ref{mouse}, bottom), particularly in genes with high expression in the SVZ tissue relative to the others. On the contrary, using the \LCLR\ transformation with the $K$-means algorithm tends to produce  tight clusters on the edges and at the vertices of the simplex (i.e., where expression in one or more of the tissues is close to zero), and a smaller number of diffuse clusters in the center of the simplex. This suggests that untransformed and CLR-transformed RNA-seq profiles tend to facilitate the identification of fine differences among nondifferentially expressed genes but do not tend to separate out genes with distinct profiles; on the other hand, the \LCLR\ transformation leads to a large diffuse grouping of nondifferentially expressed genes and the identification of several smaller, highly specific clusters. We also note that the mouse data contained a single outlier gene which exhibited aberrant expression (i.e., very strong expression in a single replicate in each tissue); only the \LCLR\ transformation grouped this gene alone in its own cluster.   

\begin{figure}[th!]\centering
\includegraphics[height=8.25cm]{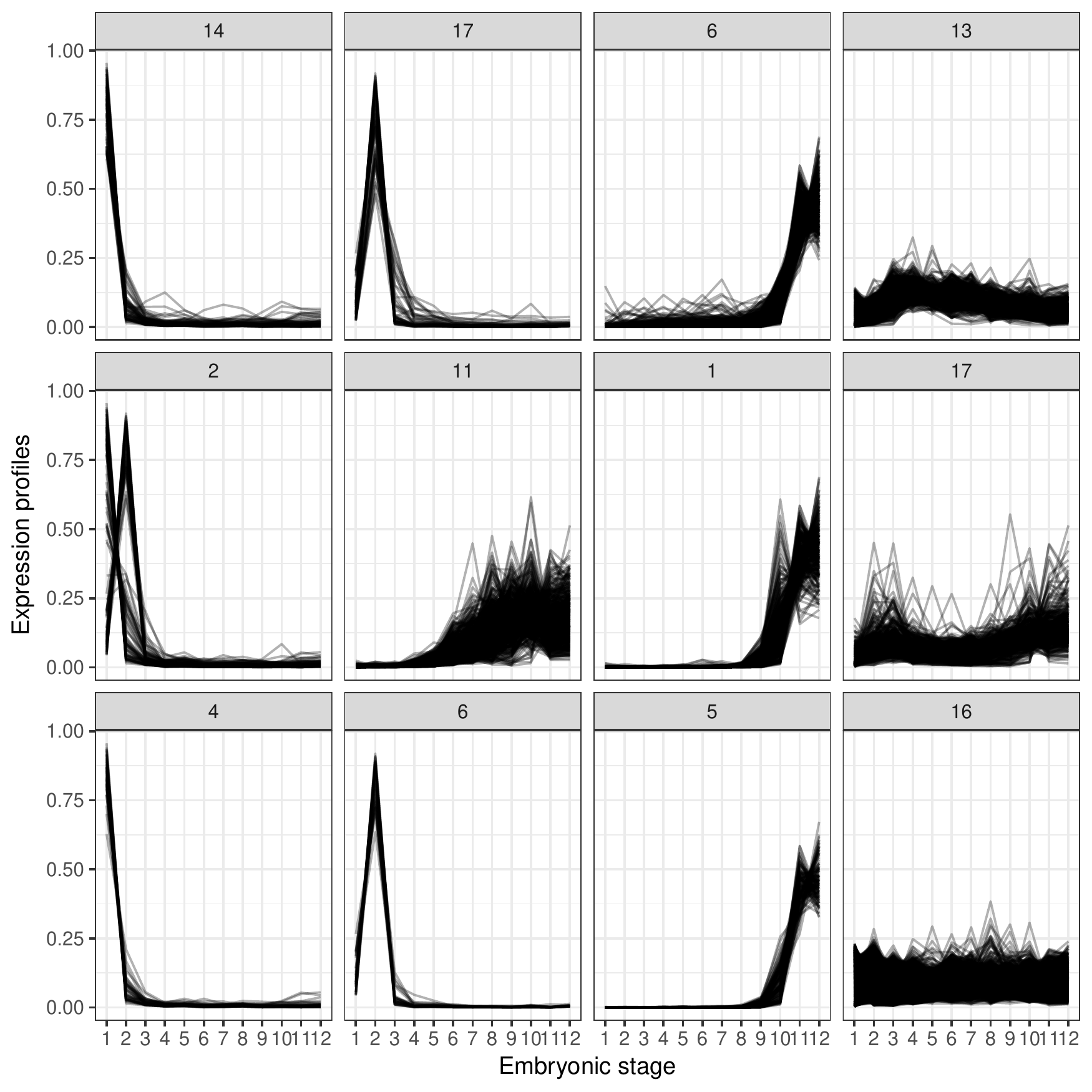}
\caption{
Examples of clusters of genes from the fly RNA-seq data: per-cluster normalized expression profiles of selected clusters obtained with the $K$-means algorithm with no transformation (top), and with the CLR (middle) and \LCLR\  (bottom) transformations are shown. Analogous clusters for each approach have been vertically aligned.\label{dro}}
\end{figure}

In the case of the fly data, the distinction between the CLR and \LCLR\ transformations is even more marked than for the mouse data. The CLR-transformed profiles with very strong relative expression in either the first or the second time point are incorrectly grouped together (Figure \ref{dro} middle, Cluster 2); these two sets of highly distinct profiles are correctly separated into two clusters in the case of no transformation (Figure \ref{dro} top, Clusters 14 and 17) or the \LCLR\ transformation (Figure \ref{dro} bottom, Clusters 4 and 6). This behavior is in fact due to the definition of the CLR transformation (see Equation~\ref{eqn:clrdef}). Indeed, as illustrated in Figure \ref{figlogclr}, the CLR transformation for coordinates close to $0$ is highly sensitive to small fluctuations; for the $K$-means algorithm, this tends to have the  effect of grouping together profiles with several zeros in common rather than those with a single non-zero coordinate in common (see for example Cluster 1 in Figure \ref{dro}, middle). In applications where it is of interest to group together profiles with common strong distinct expression, the use of untransformed or \LCLR-transformed data with $K$-means thus appears to be a more coherent choice.

\begin{figure}[th!]
\includegraphics[scale=0.28]{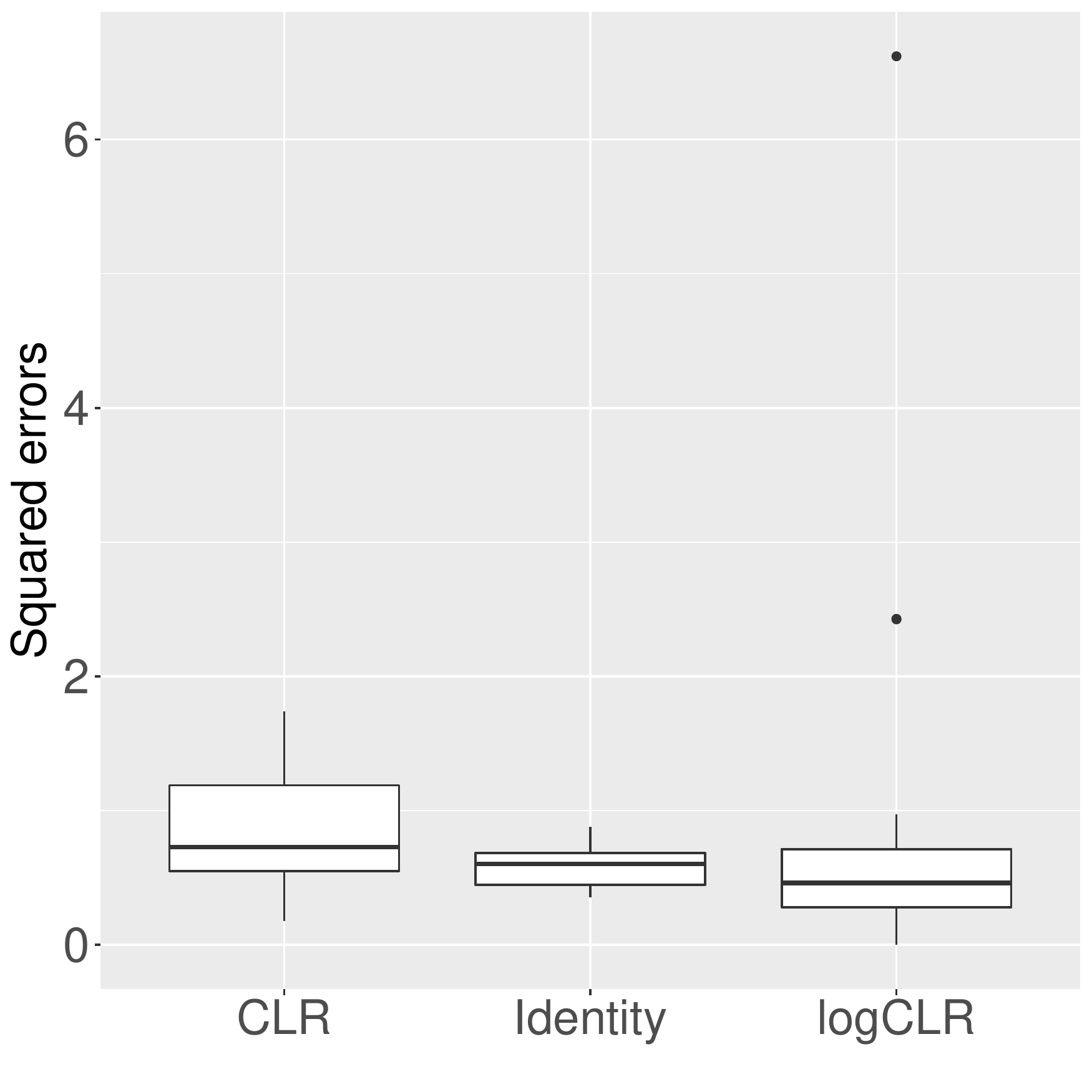}
\includegraphics[scale=0.28]{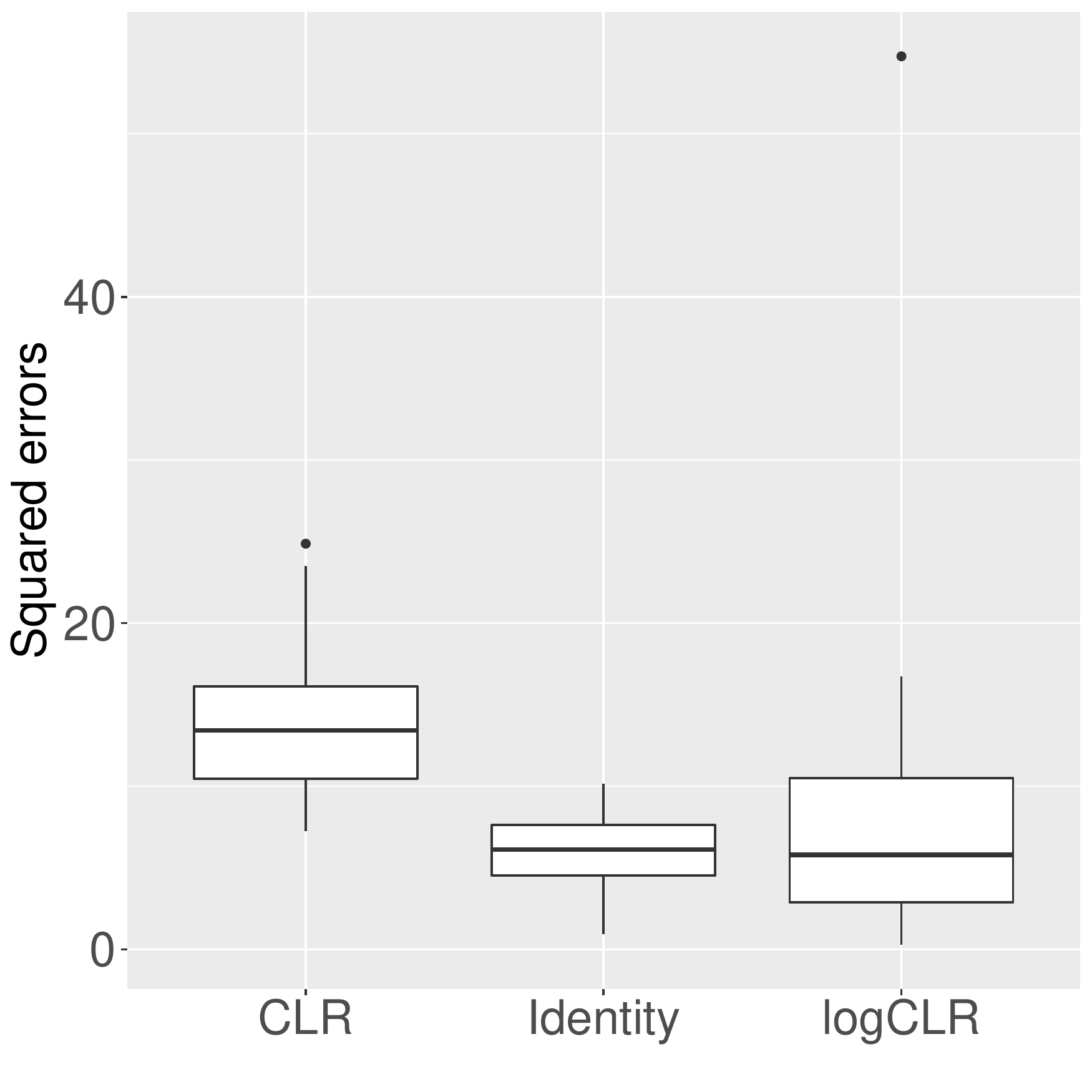}
\caption{Distributions of per-cluster squared errors for untransformed, CLR-, and \LCLR-transformed profiles for the mouse (left) and fly (right) RNA-seq data.\label{se}}
\end{figure}

To focus on a more quantitative measure of the differences among each method, we also calculate the per-cluster squared errors for each dataset using the Euclidean distance. More precisely, for each set of clusters $\mathcal{C}^{(K)}$, we consider the squared errors defined for all $k=1,\ldots,K$ by
\[
\text{SE}_{h}(k):= \sum_{i \in C_{k}} \left\| X_{i} - \mu_{k} \right\|_{2}^{2},
\]
where $\mu_{k}=\frac{1}{\vert C_{k} \vert}\sum_{i \in C_{k}} X_{i}$. As shown in Figure \ref{se}, the main difference between the \LCLR\  transformation and the other approaches is that it tends to produce a larger number of sharp clusters close to the edges or vertices of the simplex, as well as a small number of diffuse clusters situated at the center of the simplex (corresponding to the outlier points in Figure \ref{se}). As an example, Clusters $\lbrace 8,11\rbrace$ obtained from the \LCLR-transformed mouse neocortex RNA-seq data, which corresponds to a large group of nondifferentially expressed genes, is largely made up of Clusters $\lbrace 2,6,10,12,16,19,22\rbrace$ or $\lbrace 1,3,13,19\rbrace$ obtained with untransformed or CLR-transformed data, respectively. The conditional probabilities of cluster membership calculated for each approach (see the Appendix for the boxplots of conditional probabilities) indicate that these clusters tend to be made up of genes with fairly low conditional probabilities. 

\section{Bicycle sharing data results}\label{bss}

As in the previous section, we apply the $K$-means algorithm with Euclidean distance to untransformed, CLR-, and \LCLR-transformed profiles from the Velib' bicycle sharing data described in Section \ref{pres} for $K=2,\ldots,40$ clusters.  Using the nonasymptotic penalized criterion described in Section \ref{secrit}, we select a clustering model for each of the three transformations ($\widehat{K}=11,11,13$, respectively). Diagnostic plots for the calibration of the selection criterion and full profile plots for each transformation are included in the Appendix.

\begin{figure}[th!]\centering
\includegraphics[height=7.5cm]{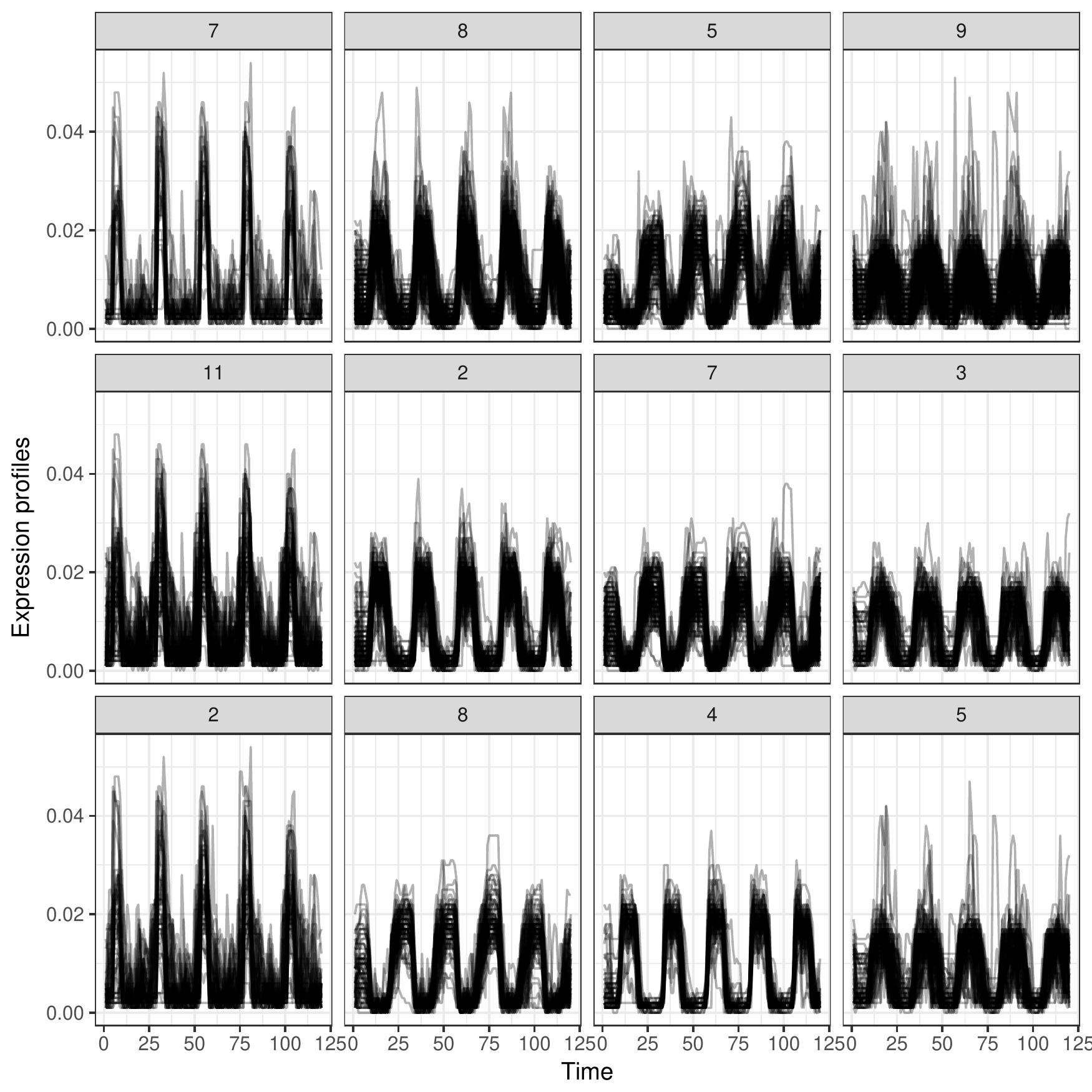}
\caption{Examples of weekly clusters of rental stations from the Velib' bicycle station data: available occupancy profiles of selected clusters obtained with the $K$-means algorithm with no transformation (top), and the CLR (middle) and \LCLR\  (bottom) transformations are shown. Analogous clusters for each approach have been vertically aligned. \label{gross}}
\end{figure}

In the case of the bicycle sharing data, the use of the CLR and \LCLR\  transformations leads to the identification of tighter clusters than for untransformed data (Figure \ref{gross}).  As was the case for RNA-seq data, the \LCLR\ transformation in particular leads to a small number of very tight clusters of distinct profiles (i.e., those near the edges or vertices of the simplex). These distinct profiles appear to largely correspond to rental stations with highly periodic rental patterns during the course of the week; for example, Cluster 4 (Figure \ref{gross}) is made up of rental stations with close to zero available docks overnight on all weeknights.

We note that the clustering partitions obtained with the CLR and \LCLR\  transformations are quite similar; the main differences between the two rely on the fact that the use of the \LCLR\  leads to the choice of a larger number of clusters with slightly tighter profiles. For clarity, in the following we thus focus on the results obtained with the CLR transformation. 

\begin{figure}[th!]\centering
\includegraphics[height=8cm]{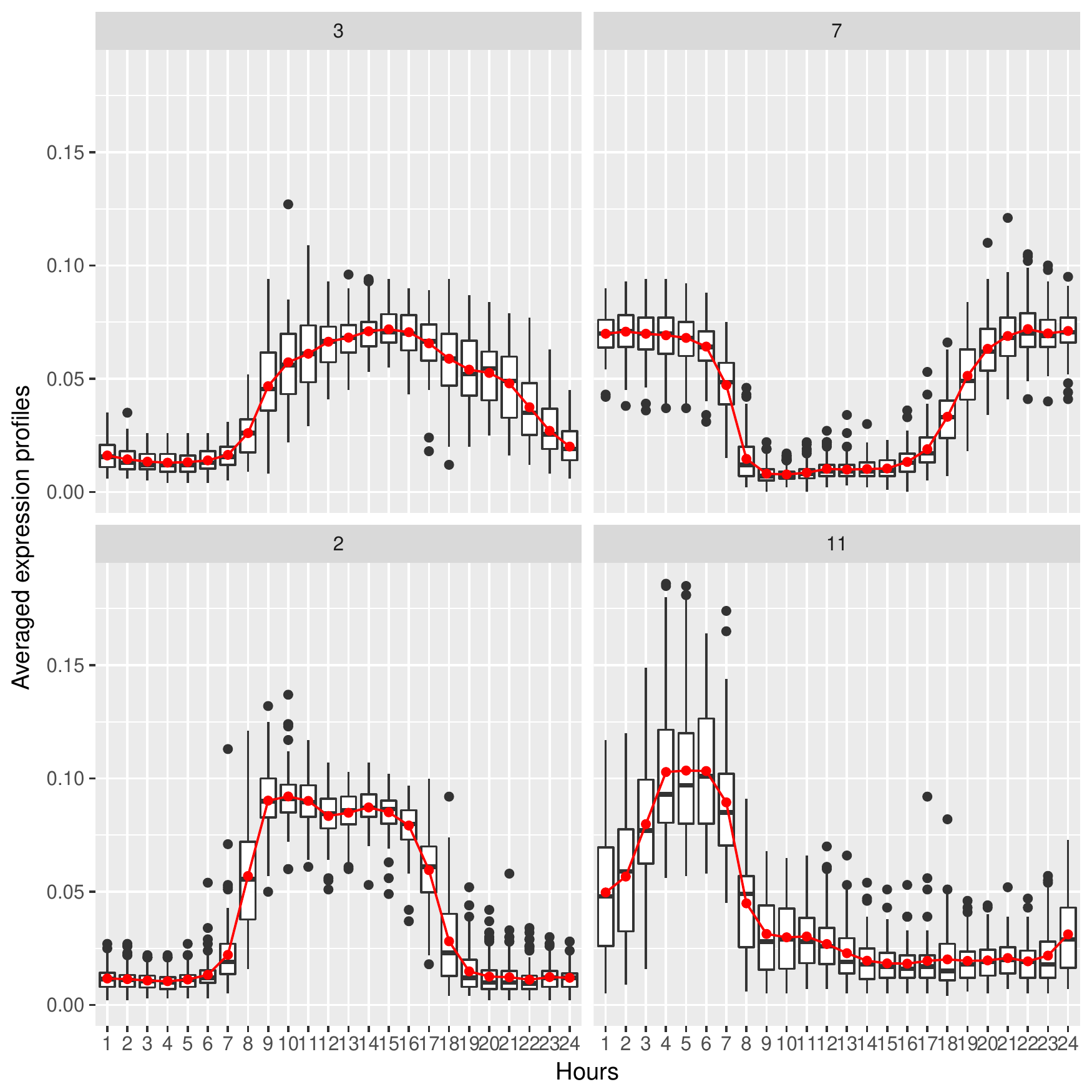}
\caption{
Examples of daily clusters of rental stations from the Velib' bicycle station data: available occupancy profiles of selected clusters obtained with the $K$-means algorithm CLR transformation. Connected red lines correspond to the mean profiles for each time point. 
 \label{velib_clr}}
\end{figure}

\begin{figure}[th!]\centering
\includegraphics[height=8cm]{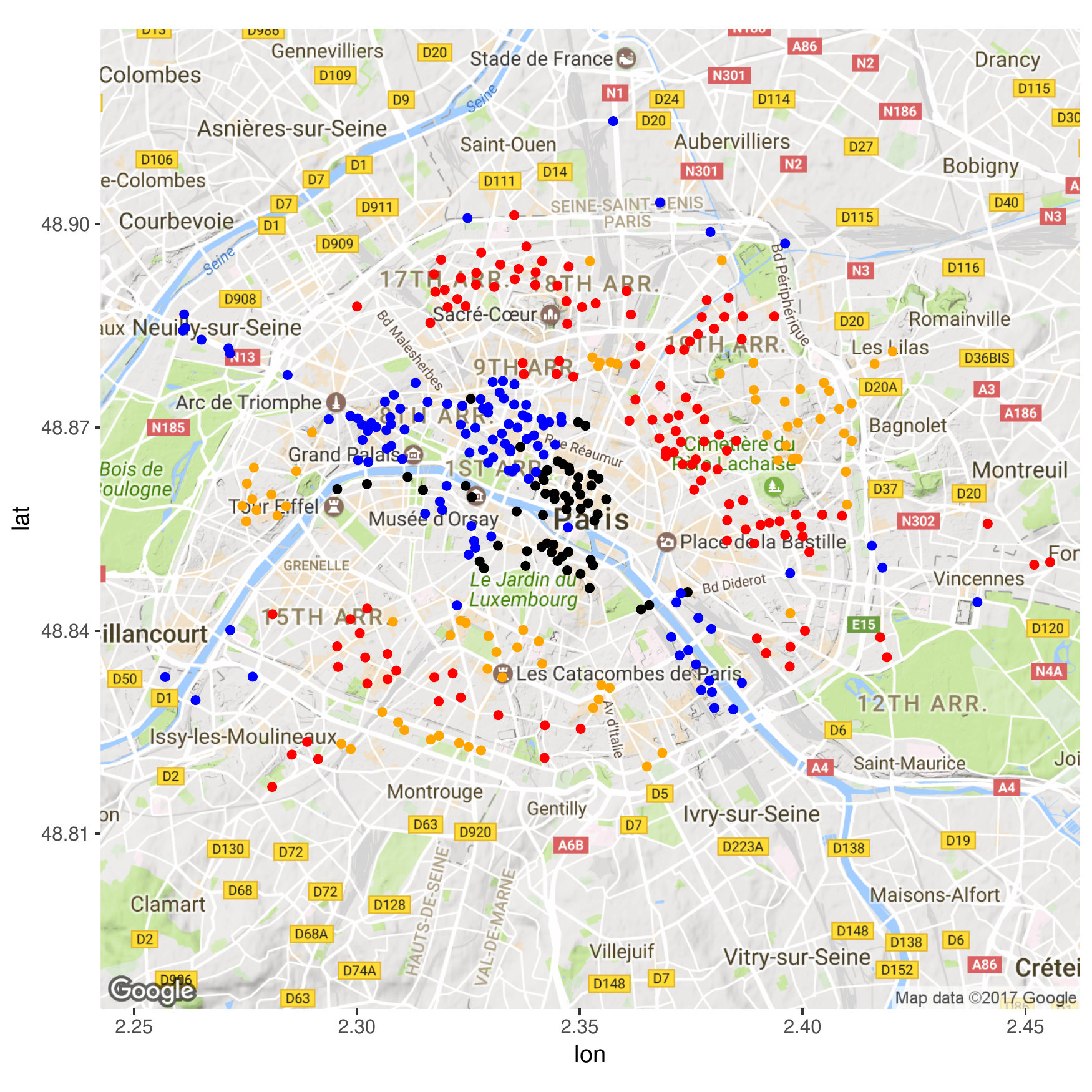}
\caption{Geographic location in Paris of selected rental station clusters from the Velib' bicycle station occupancy data obtained with the K-means algorithm and CLR transformation. 
\label{map}}
\end{figure}

As shown in Figure \ref{velib_clr}, the daily cluster profiles (i.e., after summing occupancy profiles for each time point across the five days) reveal an interesting coherence with the location of the Velib' rental stations; this is all the more striking given that the spatial coordinates of the stations are not included in the model. For example, Cluster 3 (Figure \ref{velib_clr}, top left) corresponds to stations located in the historical city center of Paris (in black in Figure \ref{map}). This area of Paris is a shopping and tourist-oriented district with many restaurants and bars, which explains the fact that users appear to arrive over the course of the morning and to have prolonged departures over the course of the evening. Cluster 7 (Figure \ref{velib_clr}, top right) corresponds to stations serving users that leave for work in the morning and return home at the end of the workday. These stations (in red in Figure \ref{map}) indeed correspond to largely residential areas of Paris. On the other hand, Cluster 2 (Figure \ref{velib_clr}, bottom left) appears to contain rental stations for professionals arriving at business offices, as users tend to arrive between $8$am and $10$am and leave between $6$pm to $8$pm. This is again confirmed by the map, where these stations (in blue) are located in main business districts such as Bercy, La D\'efense, and the $8$th Arrondissement. Finally, Cluster $11$ (Figure \ref{velib_clr}, bottom right) appears to represent a fairly atypical group of rental stations, as they fill up early in the morning, between $1$am to $5$am; in addition, their corresponding locations (in orange in Figure \ref{map}) are largely residential (especially in the $19$th and $20$th Arrondissements) and do not necessarily match up to areas of active night life. However, most of these stations belong to the special category of V+ stations located at higher elevations; as such, users tend to rent bicycles in the morning (on the downhill route) but do not return home by bicycle, thus necessitating an overnight redistribution of bicycles. 

\section{Conclusion and recommendations}\label{conclusion}
In this work, we have described a strategy for clustering compositional data with the $K$-means algorithm and several adaptated transformations. The choice of the appropriate transformation in practice depends strongly on the type of cluster profiles that are of interest for a given context; for example, is it pertinent to separate highly-specific profiles into dedicated clusters (as appears to be the case for the RNA-seq data)? Our recommendations for the choice of compositional data transformation can be summarized into the following points: 1) if a balanced clustering is desired (i.e., clusters of roughly equal size), the $K$-means algorithm should be applied on untransformed profiles; 2) in order to highlight groups of highly-specific profiles, the \LCLR\ transformation should instead be used; and 3) if small fluctuations of profile coordinates close to 0 are of primary interest, the CLR transformation may be preferred over the \LCLR. Recall that in the two applications studied here, the \LCLR\ transformation was the only option that yielded satisfactory results in all cases.

Since the data studied here originated as count tables that were transformed into compositional tables, it is worth noting that a $K$-means algorithm with the $\chi^2$ distance could also have been applied on the original counts. This distance is based on a transformation that relies on the multiplication of columns by constants which are, in the case of the applications studied here, very close to one another (in the case of the RNA-seq data, this is due to the normalization of library sizes). As such, results obtained for the $K$-means algorithm with Euclidean distance on CLR- or \LCLR-transformed profiles are very similar to those obtained for the $K$-means algorithm with $\chi^2$ distance on the original counts.  

Finally, several extensions of this work could be considered in future work. For instance, rather than selecting a single number of clusters for a given dataset, it could in some cases be preferable to instead construct a hierarchy of clustering for varying $K$. To do this, one possibility would be to partition the data into a moderate number of clusters, for example using the strategy described in this work, prior to aggregating clusters according to an appropriate criterion.  Another potential extension is the integration of relevant external data in the clustering approach, such as location (for the Velib' data) or membership in functional pathways (for the RNA-seq data).

\def\cprime{$'$}

\newpage
\begin{appendix}
\section{Results for the mouse RNA-seq data}
In the following, we visualize the clusters of genes identified from the mouse neocortex RNA-seq data after applying the $K$-means algorithm with the identity, CLR, and \LCLR\ transformations, as well as plots of the conditional probabilities calculated for each. 
We also present an example of the slope heuristic plot (obtained with the R package \texttt{capushe}) for calibrating of the nonasymptotic penalized criterion defined in Section 3.2 of the main manuscript.

\begin{figure}[H]
\includegraphics[width=13cm,height=13cm]{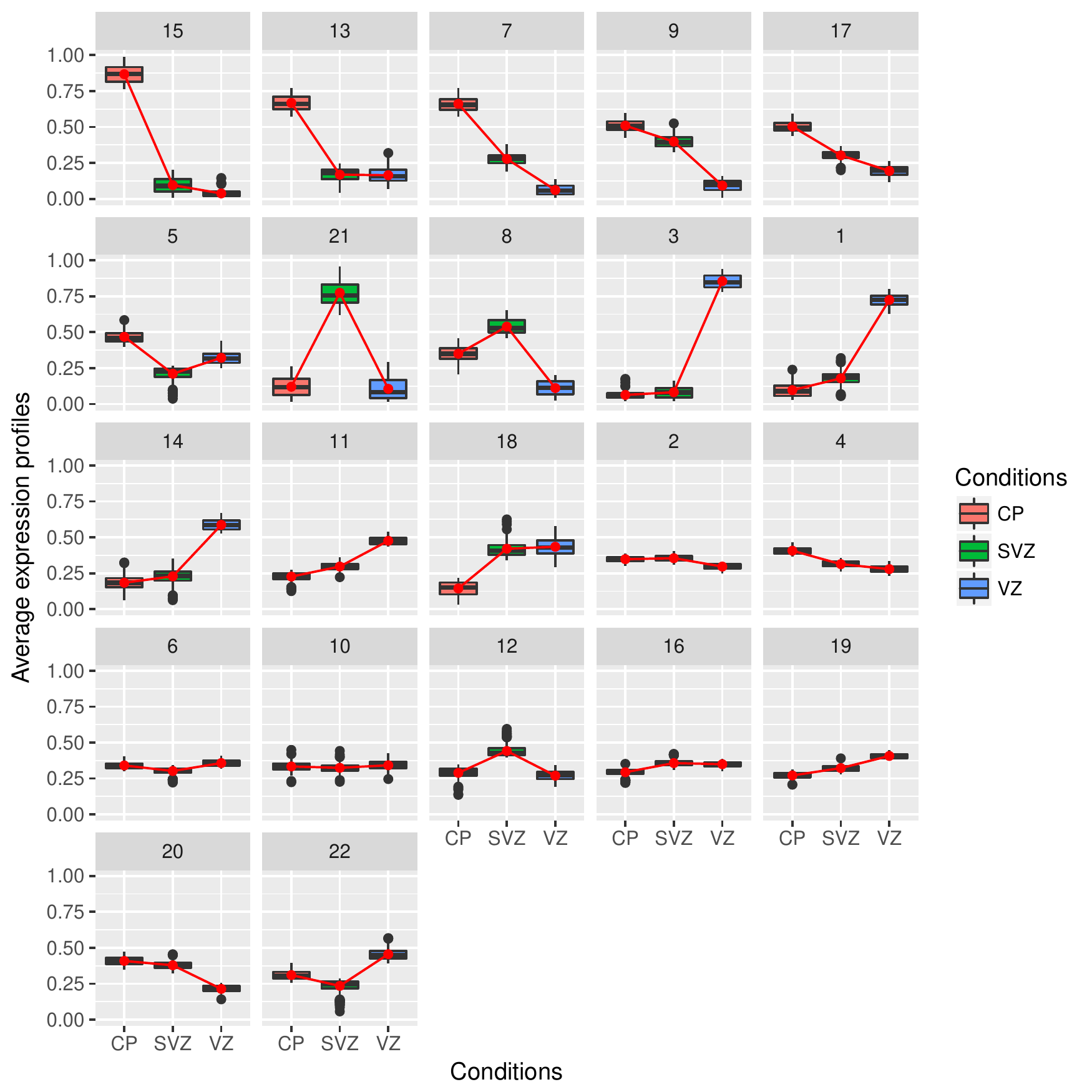}
\caption{Clusters of genes from the mouse neocortex RNA-seq data: per-cluster normalized expression profiles of selected clusters obtained with the $K$-means algorithm and no transformation. Clusters have been arranged so that those with similar average profiles are displayed next to one another. Connected red lines correspond to the mean profiles for each tissue.}
\end{figure}

\begin{figure}[H]
\includegraphics[width=15cm,height=15cm]{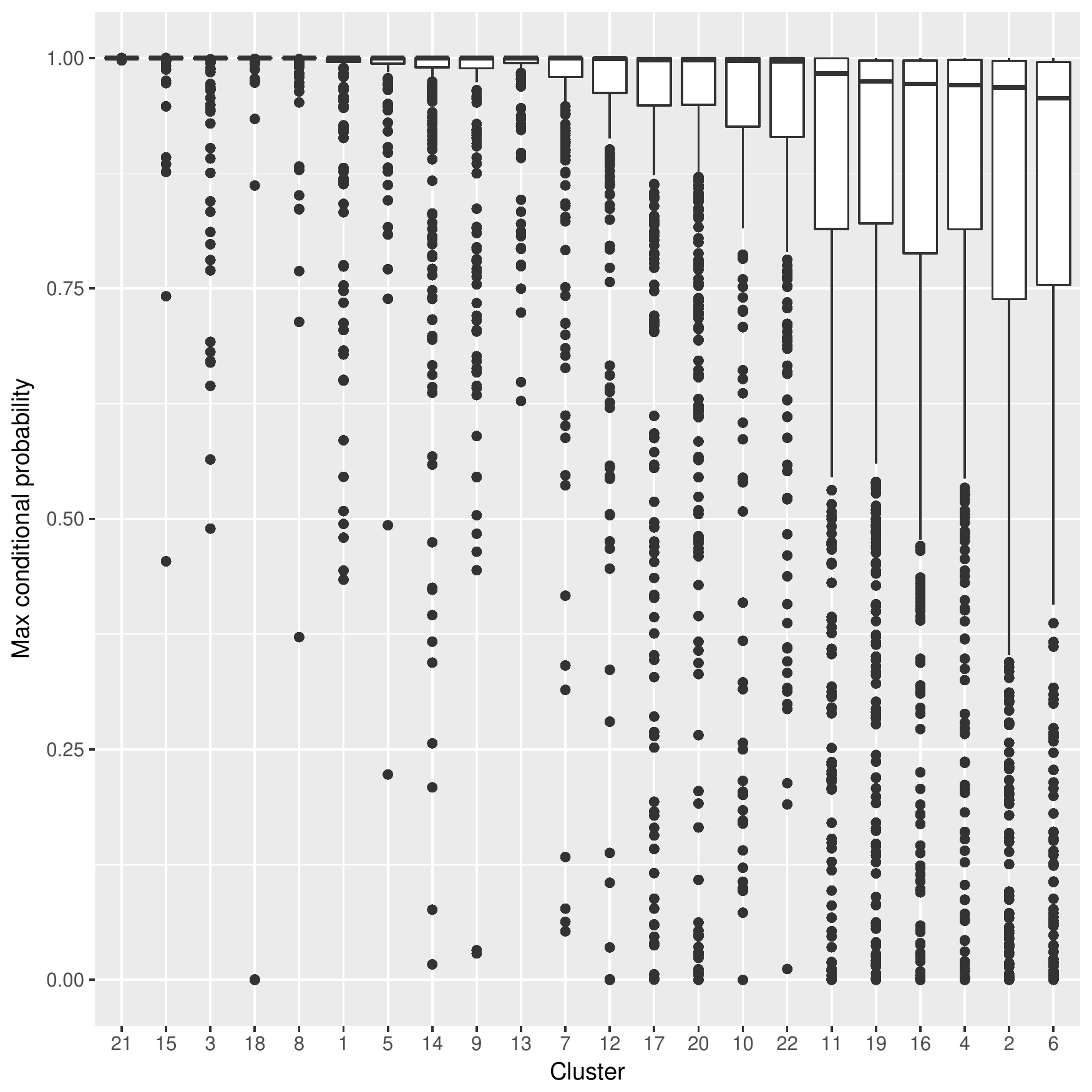}
\caption{Per-cluster maximum conditional probabilities from the mouse neocortex RNA-seq data obtained with the $K$-means algorithm and no transformation. Clusters have been  sorted by median per-cluster conditional probabilities. }
\end{figure}

\begin{figure}[H]
\includegraphics[width=15cm,height=15cm]{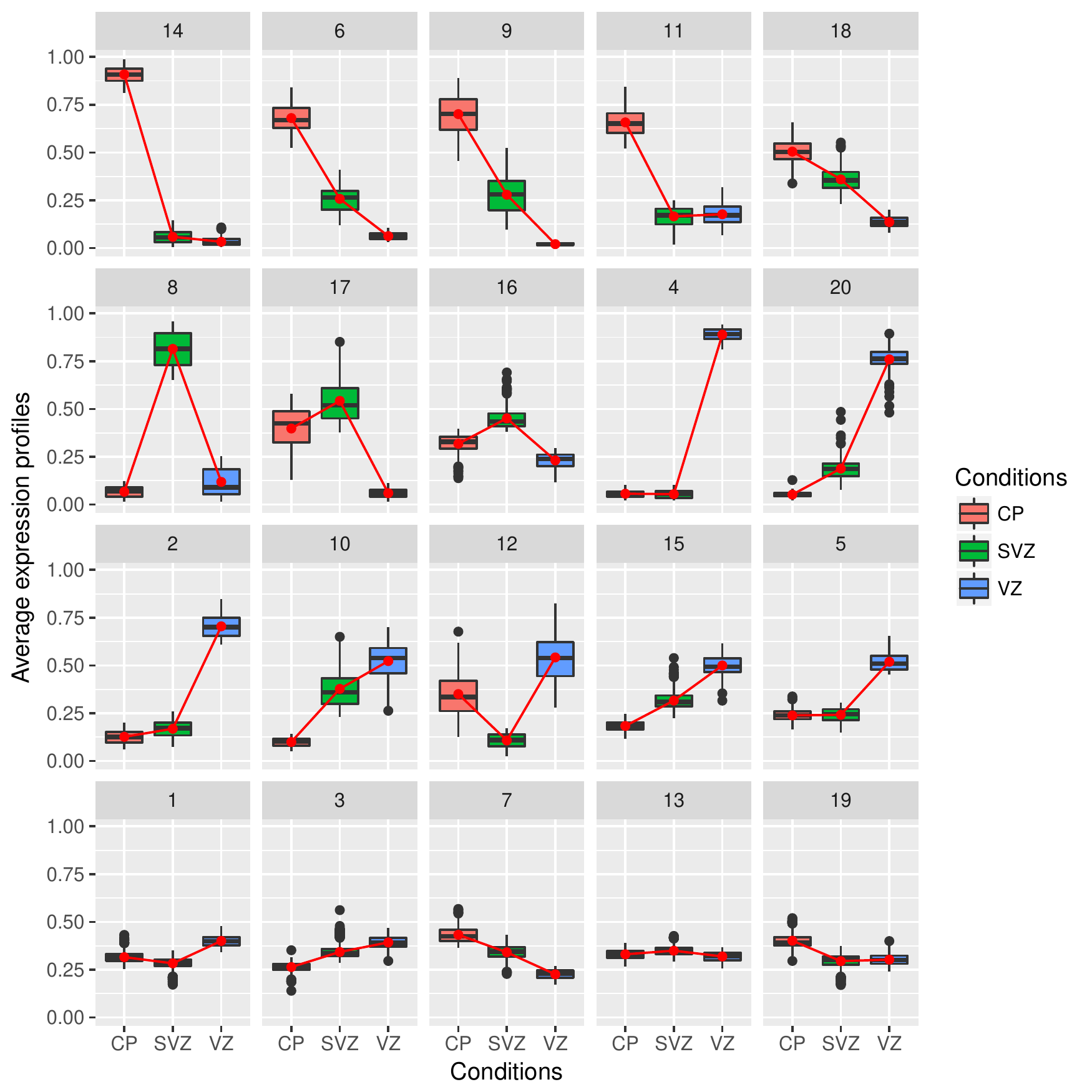}
\caption{Clusters of genes from the mouse neocortex RNA-seq data: per-cluster normalized expression profiles of selected clusters obtained with the $K$-means algorithm and CLR transformation. Clusters have been arranged so that those with similar average profiles are displayed next to one another. Connected red lines correspond to the mean profiles for each tissue.}
\end{figure}

\begin{figure}[H]
\includegraphics[width=15cm,height=15cm]{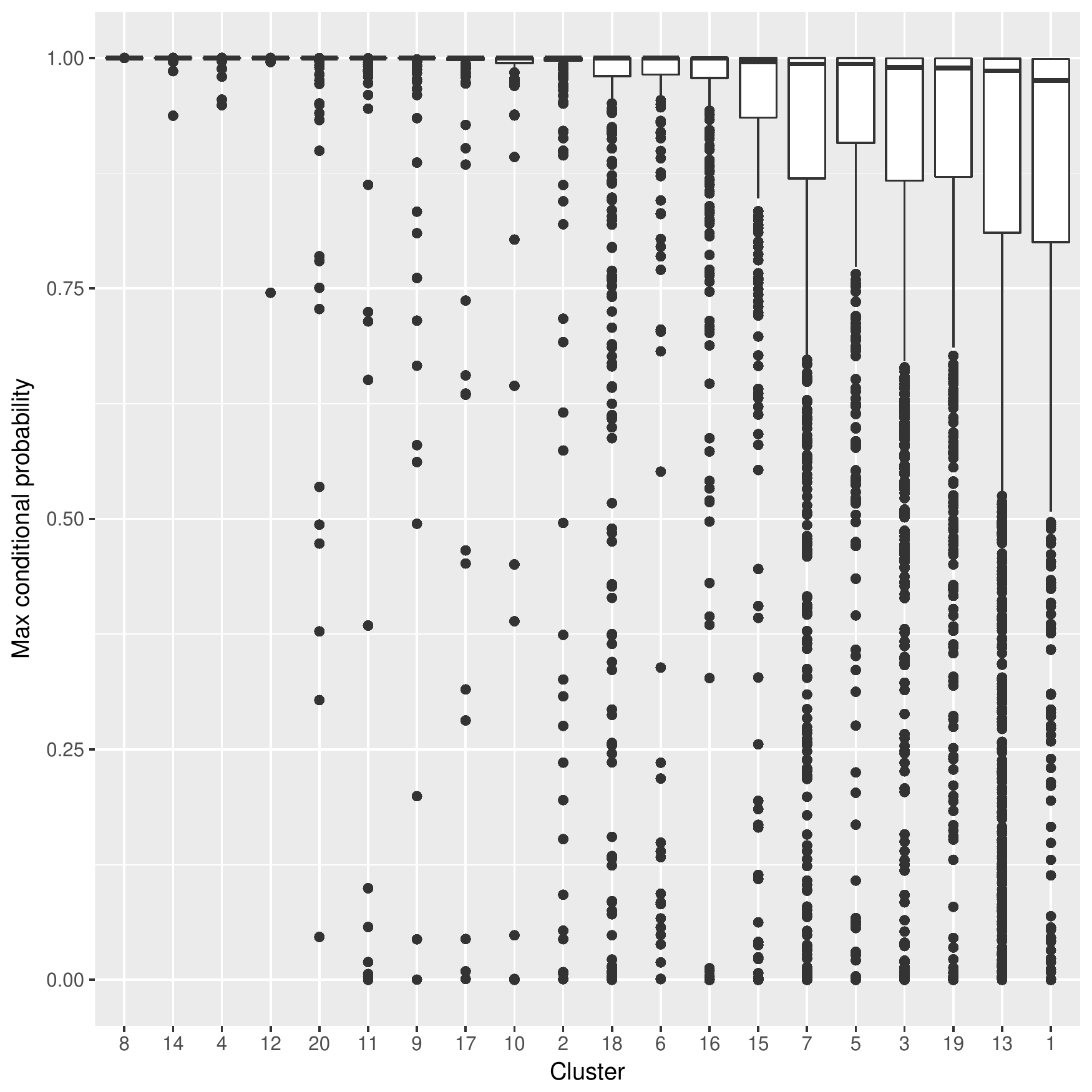}
\caption{Per-cluster maximum conditional probabilities from the mouse neocortex RNA-seq data obtained with the $K$-means algorithm and CLR transformation. Clusters have been  sorted by median per-cluster conditional probabilities.}
\end{figure}

\begin{figure}[H]
\includegraphics[width=15cm,height=15cm]{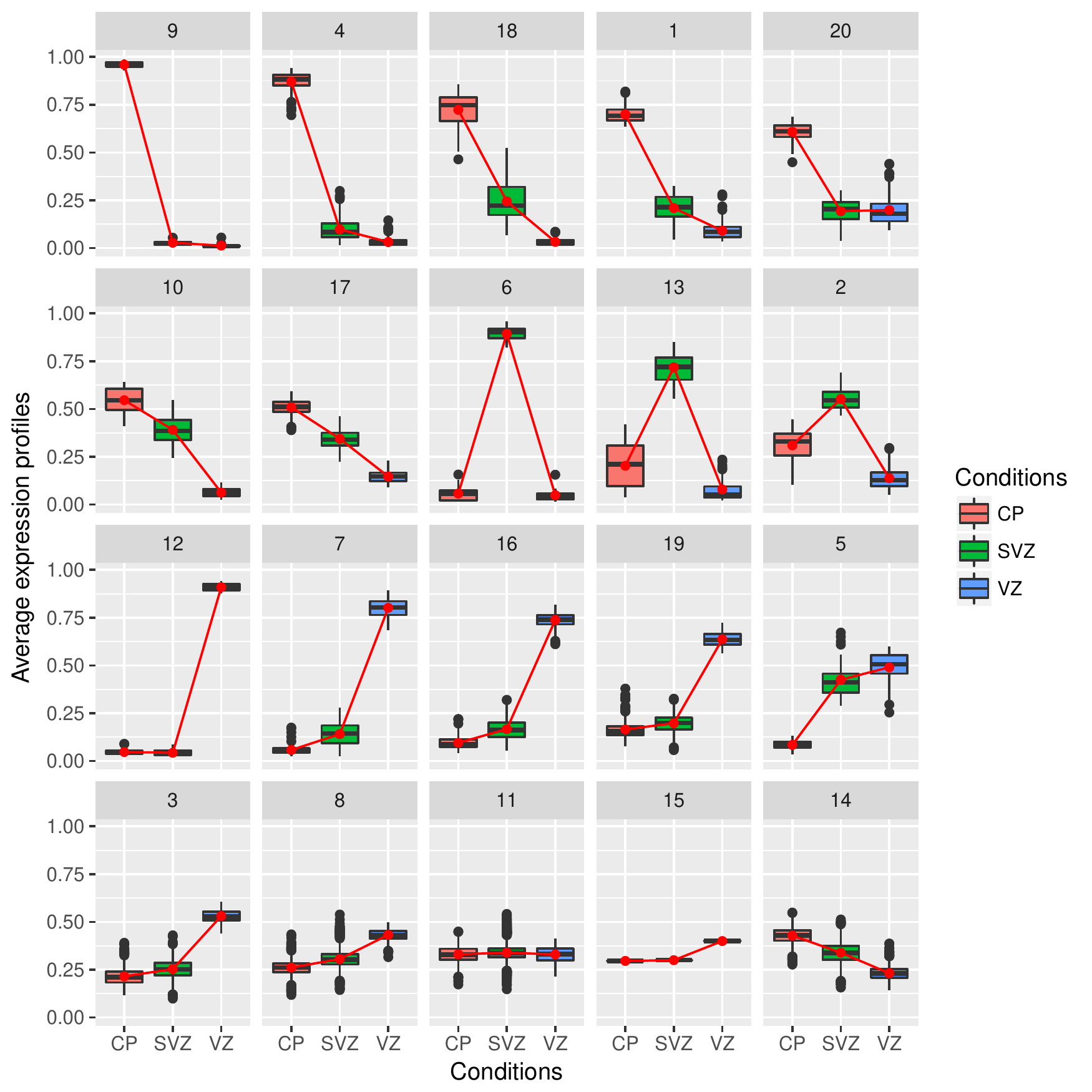}
\caption{Clusters of genes from the mouse neocortex RNA-seq data: per-cluster normalized expression profiles of selected clusters obtained with the $K$-means algorithm and \LCLR\ transformation. Clusters have been arranged so that those with similar average profiles are displayed next to one another. Connected red lines correspond to the mean profiles for each tissue.}
\end{figure}

\begin{figure}[H]
\includegraphics[width=15cm,height=15cm]{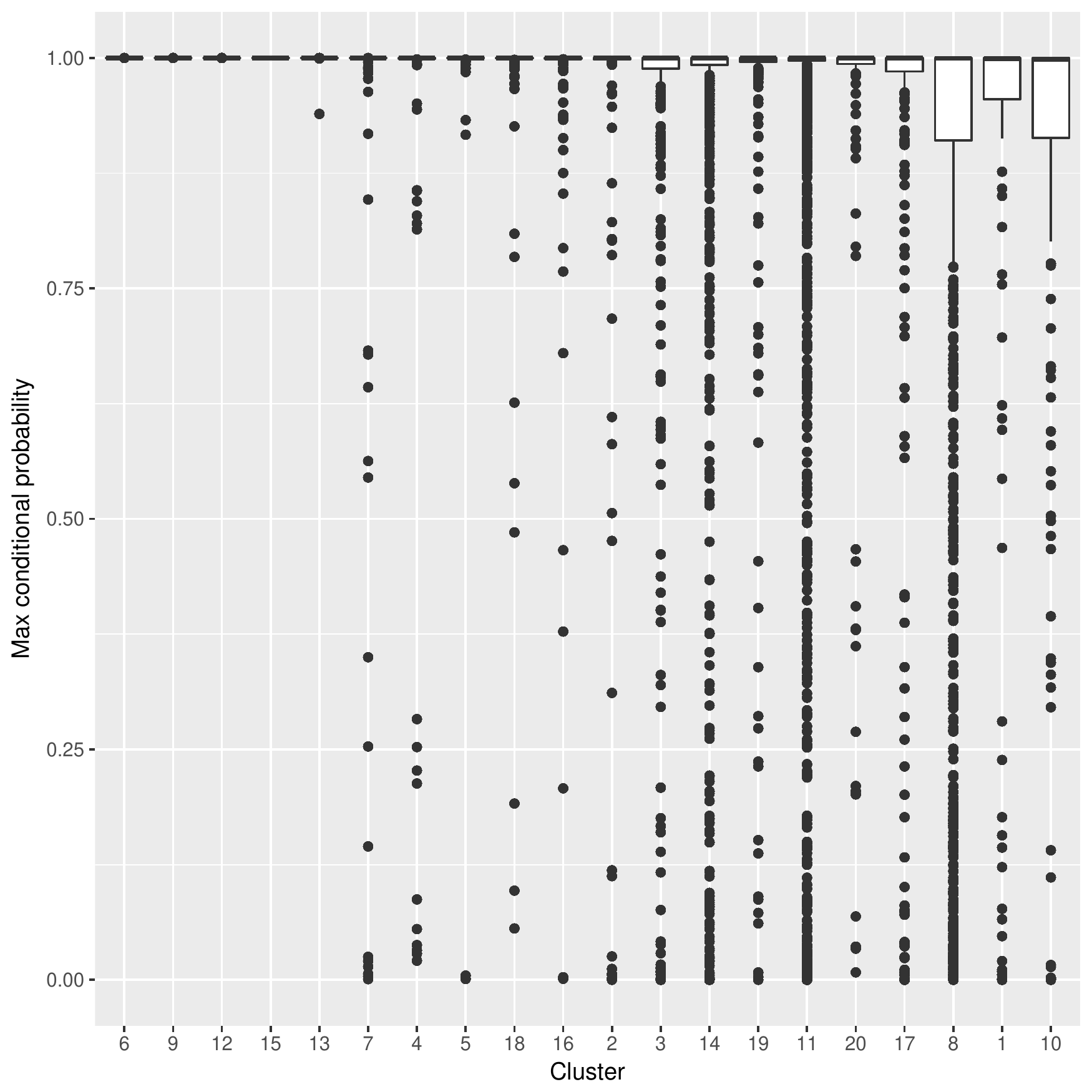}
\caption{Per-cluster maximum conditional probabilities from the mouse neocortex RNA-seq data obtained with the $K$-means algorithm and \LCLR\ transformation. Clusters have been  sorted by median per-cluster conditional probabilities.}
\end{figure}

\begin{figure}[H]
\includegraphics[width=15cm,height=15cm]{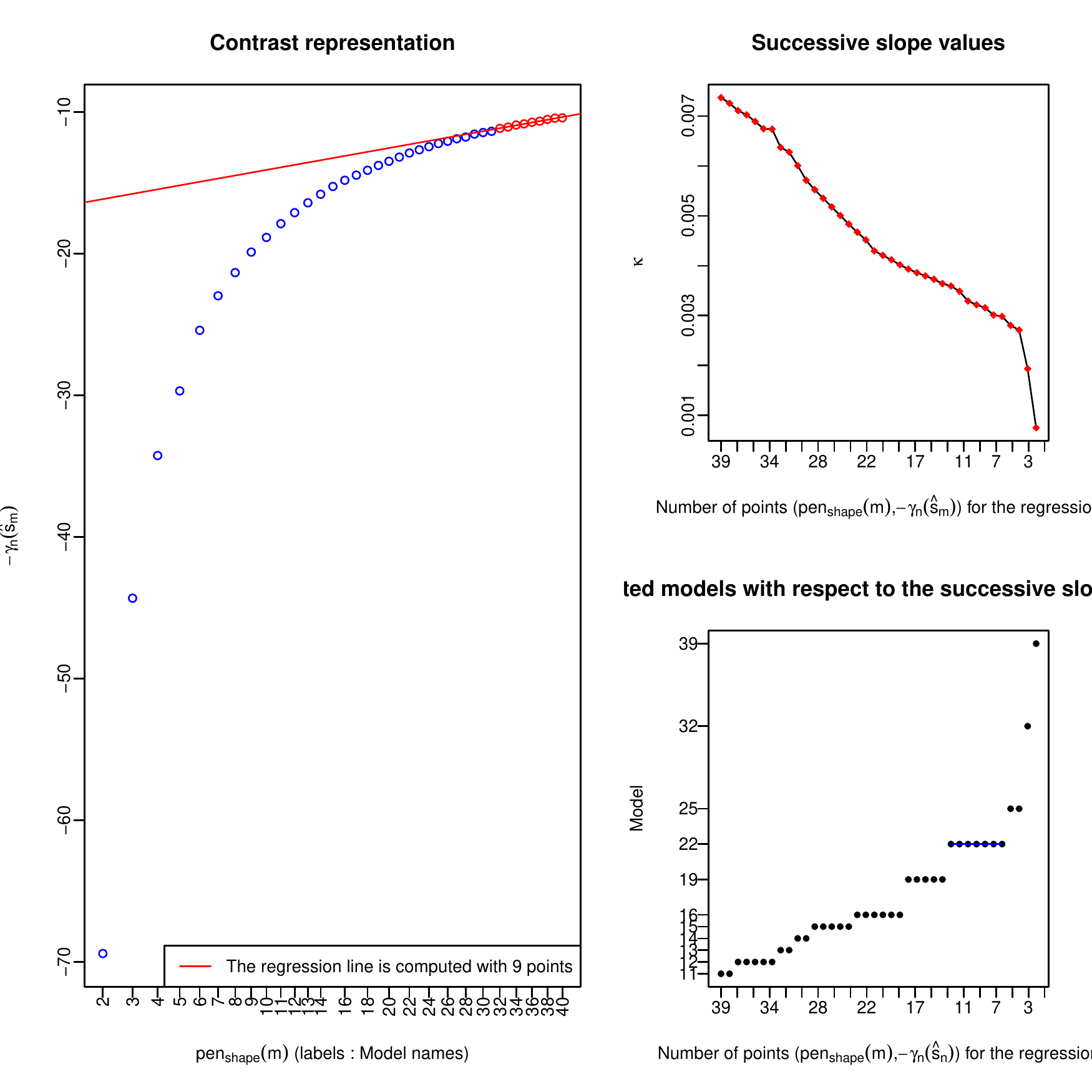}
\caption{Example of slope heuristics for clustering the mouse neocortex RNA-seq data with the $K$-means algorithm and no transformation.}
\end{figure}


\section{Results for the fly embryonic RNA-seq data}
In the following, we visualize the clusters of genes identified from the fly embryonic  RNA-seq data after applying the $K$-means algorithm with the identity, CLR, and \LCLR\ transformations.

\begin{figure}[H]
\includegraphics[width=15cm,height=15cm]{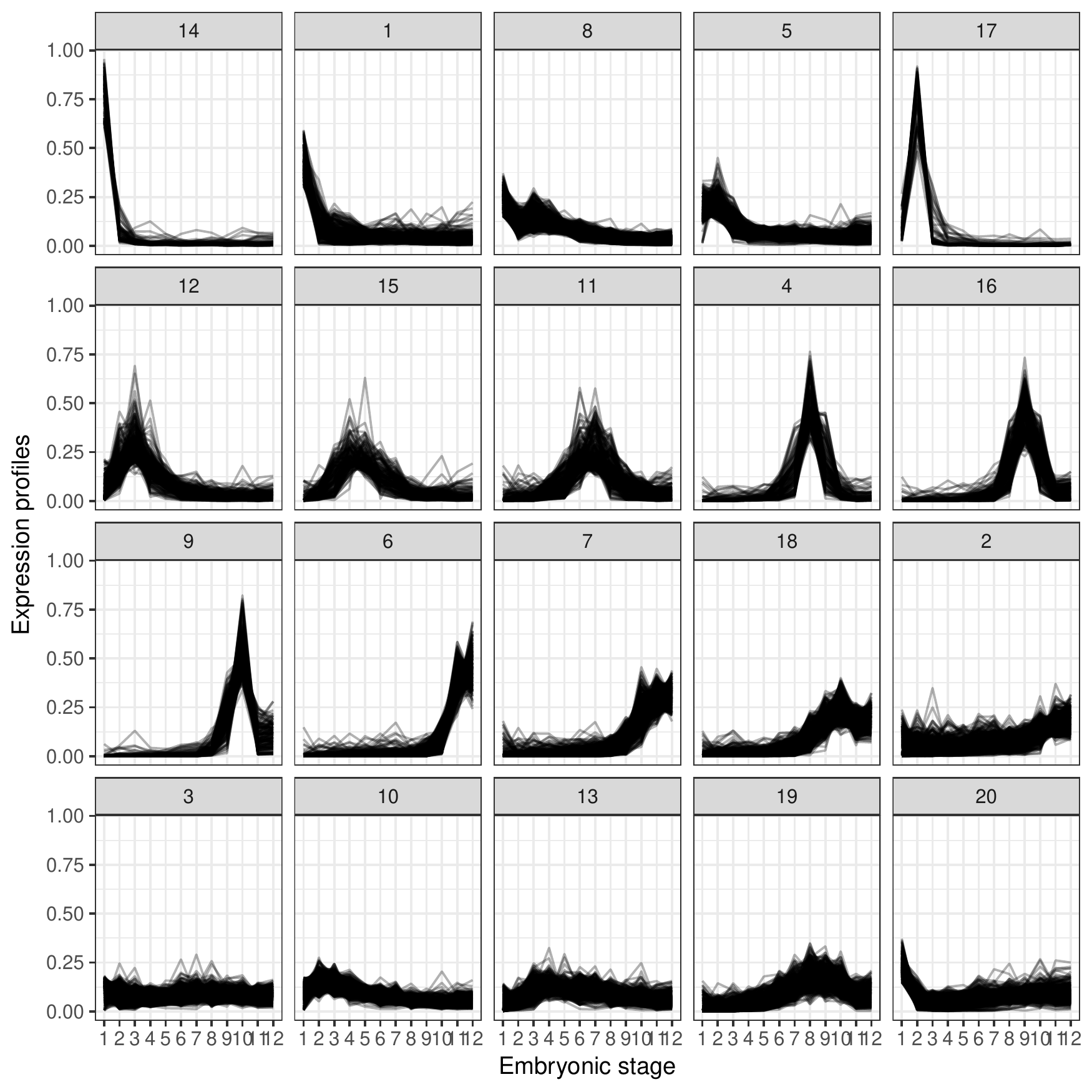}
\caption{Clusters of genes from the fly embryonic RNA-seq data: per-cluster normalized expression profiles of selected clusters obtained with the $K$-means algorithm and no transformation. Clusters have been arranged so that those with similar average profiles are displayed next to one another.}
\end{figure}

\begin{figure}[H]
\includegraphics[width=15cm,height=15cm]{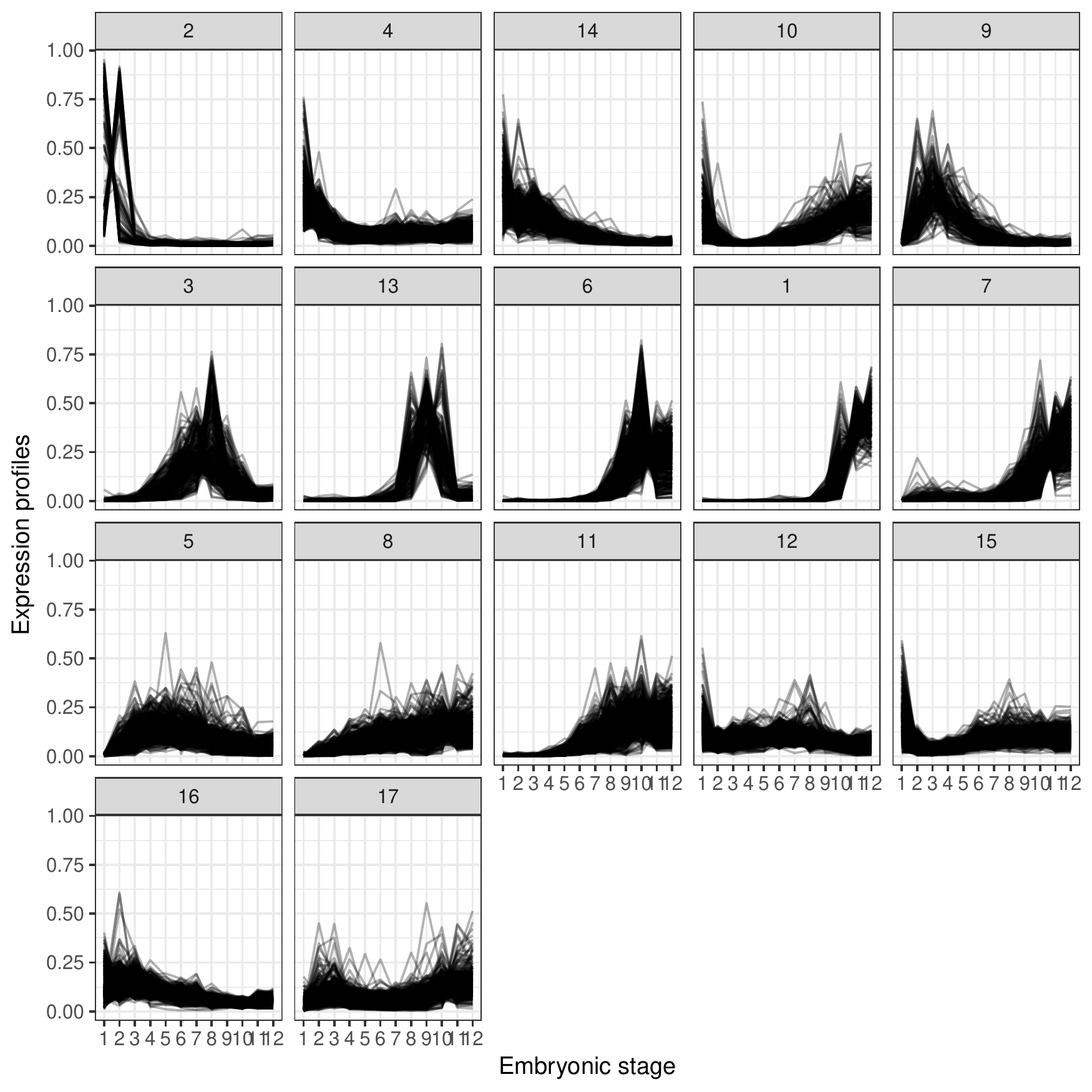}
\caption{Clusters of genes from the fly embryonic RNA-seq data: per-cluster normalized expression profiles of selected clusters obtained with the $K$-means algorithm and CLR transformation. Clusters have been arranged so that those with similar average profiles are displayed next to one another.}
\end{figure}

\begin{figure}[H]
\includegraphics[width=15cm,height=15cm]{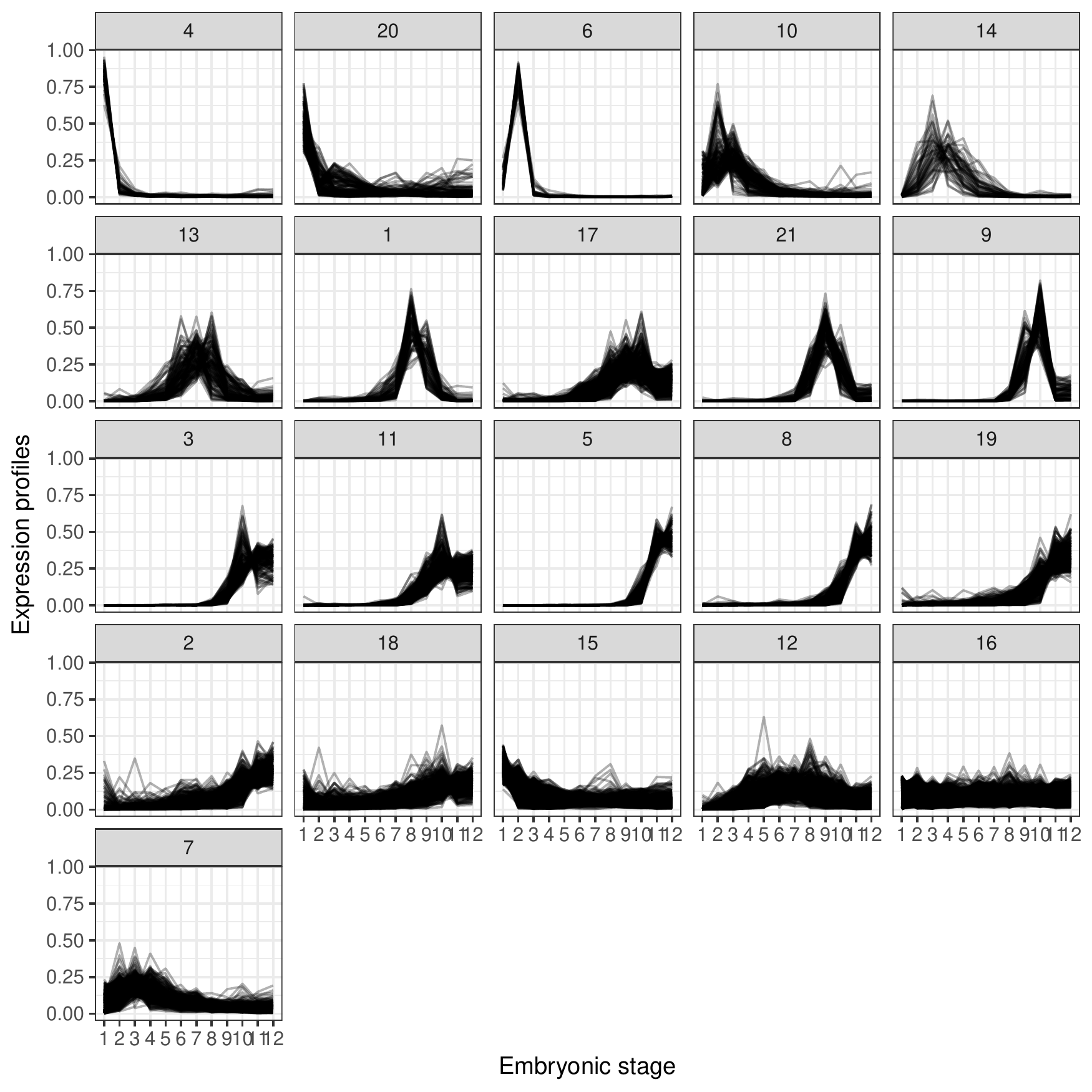}
\caption{Clusters of genes from the fly embryonic RNA-seq data: per-cluster normalized expression profiles of selected clusters obtained with the $K$-means algorithm and \LCLR\ transformation. Clusters have been arranged so that those with similar average profiles are displayed next to one another.}
\end{figure}


\section{Results for the Velib' bicycle sharing system data}

In the following, we visualize the clusters of genes identified from the Velib' bicycle sharing system data after applying the $K$-means algorithm with the identity, CLR, and \LCLR\ transformations.

\begin{figure}[H]
\includegraphics[width=15cm,height=15cm]{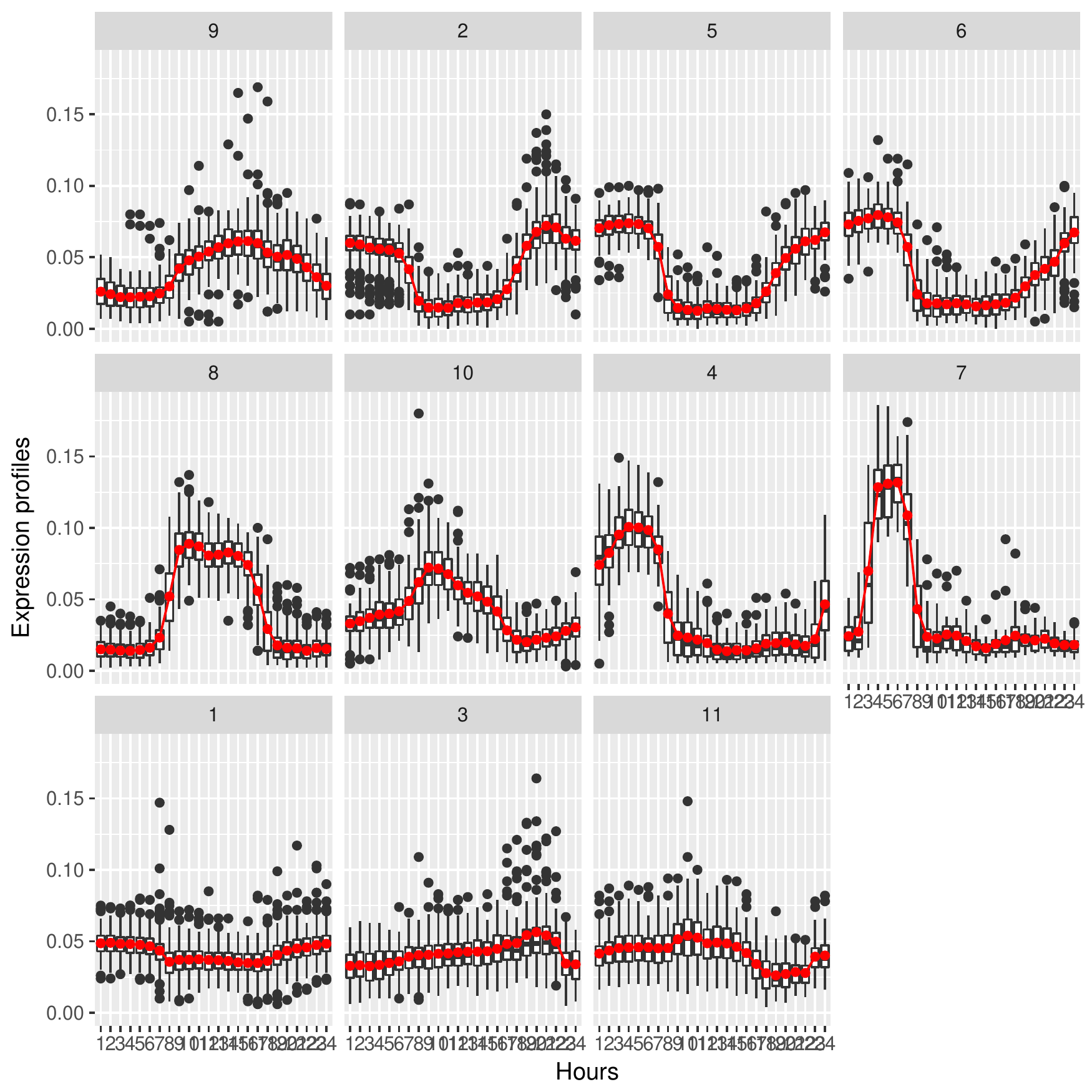}
\caption{Clusters of genes from the Velib' bicycle sharing system data: per-cluster occupancy profiles of selected clusters obtained with the $K$-means algorithm and no transformation. Clusters have been arranged so that those with similar average profiles are displayed next to one another. Connected red lines correspond to the mean profiles for each time point.}
\end{figure}

\begin{figure}[H]
\includegraphics[width=15cm,height=15cm]{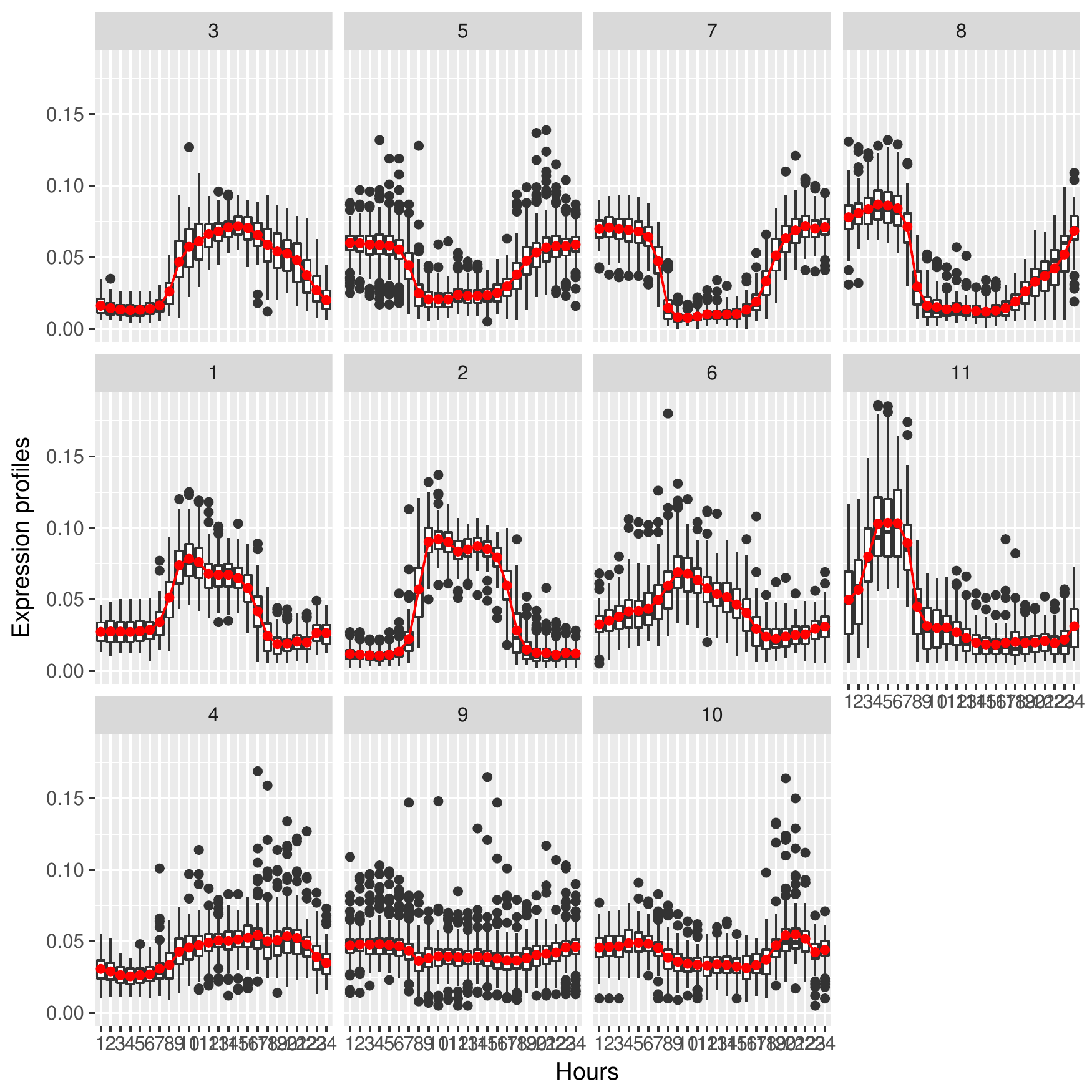}
\caption{Clusters of genes from the Velib' bicycle sharing system data: per-cluster occupancy profiles of selected clusters obtained with the $K$-means algorithm and CLR transformation. Clusters have been arranged so that those with similar average profiles are displayed next to one another. Connected red lines correspond to the mean profiles for each time point.}
\end{figure}

\begin{figure}[H]
\includegraphics[width=15cm,height=15cm]{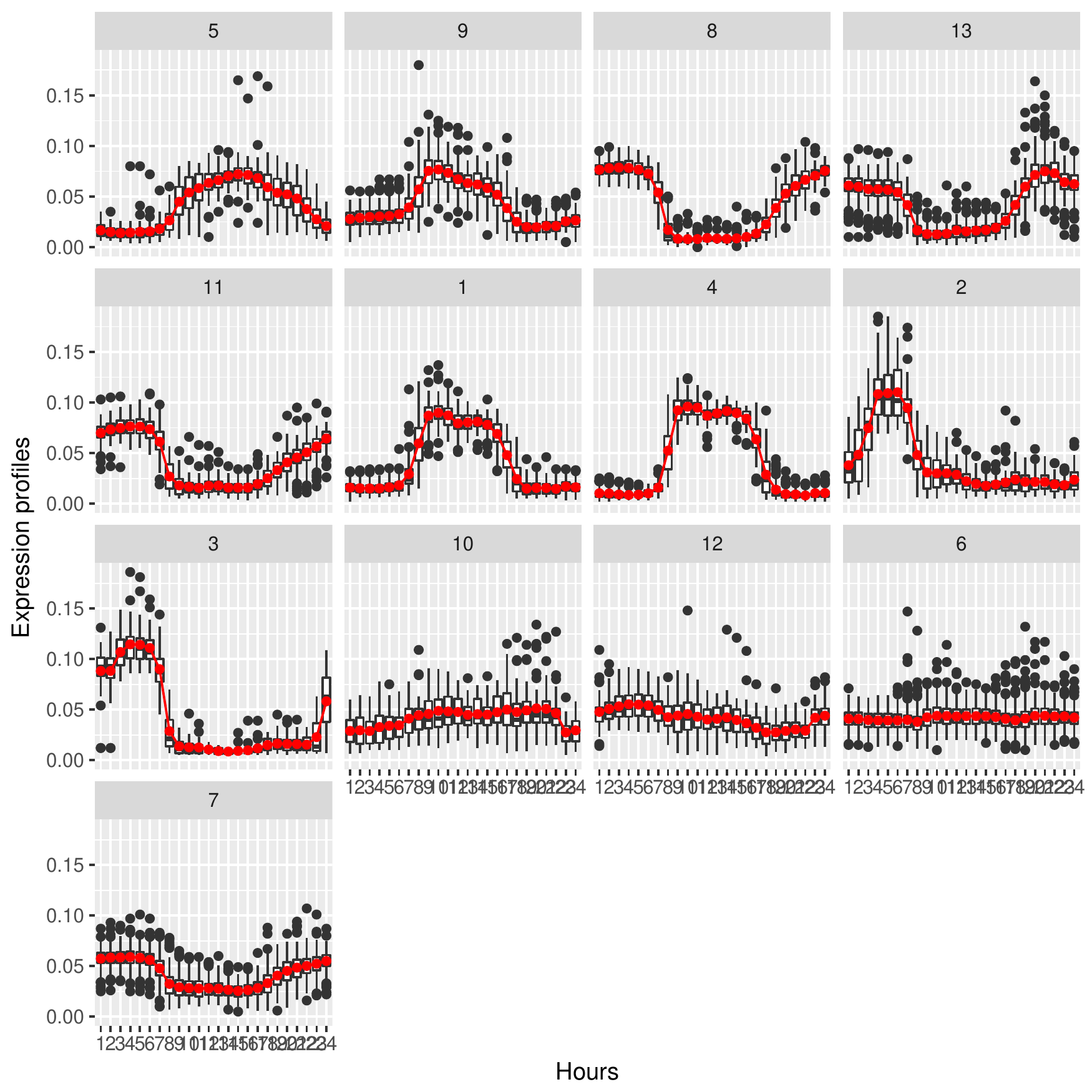}
\caption{Clusters of genes from the Velib' bicycle sharing system data: per-cluster occupancy profiles of selected clusters obtained with the $K$-means algorithm and \LCLR\ transformation. Clusters have been arranged so that those with similar average profiles are displayed next to one another. Connected red lines correspond to the mean profiles for each time point.}
\end{figure}

\end{appendix}

\end{document}